\documentclass{amsart}

\usepackage{amssymb}

\numberwithin{equation}{subsection}

\theoremstyle{plain}

\newtheorem{theorem}[equation]{Theorem}

\newtheorem{proposition}[equation]{Proposition}

\newtheorem{hypo}[equation]{Hypothesis}

\newtheorem{lemma}[equation]{Lemma}

\newtheorem{claim}[equation]{Claim}

\newtheorem{case}[equation]{Case}

\newtheorem*{teoremaa}{Theorem A}

\newtheorem*{teoremab}{Theorem B}

\newcommand{\Irr}{\operatorname{Irr}}

\newcommand{\cd}{\operatorname{cd}}

\newcommand{\dl}{\operatorname{dl}}

\newcommand{\Gal}{\operatorname{Gal}}

\newcommand{\Ker}{\operatorname{Ker}}

\newcommand{\Aut}{\operatorname{Aut}}

\newcommand{\Inn}{\operatorname{Inn}}

\newcommand{\modu}{\operatorname{mod}}

\newcommand{\Lin}{\operatorname{Lin}}

\newcommand{\GL}{\operatorname{GL}}

\newcommand{\GF}{\operatorname{GF}}

\newcommand{\SL}{\operatorname{SL}}

\newcommand{\Spl}{\operatorname{Sp}}

\theoremstyle{definition}

\newtheorem{definition}[equation]{Definition}

\begin{document}

	\title{Products of characters with few irreducible constituents}

\author{Edith Adan-Bante}

\address{University of Southern Mississippi at Gulf Coast,
Long Beach, MS, 39560}

\email{Edith.Bante@usm.edu}

\keywords{Product of characters, irreducible character,  derived length}

\thanks{This research was partially supported by grant NSF DMS 9970030.}

\subjclass{20c15}

\date{2004}

\begin{abstract}

 We 
study the solvable groups $G$ that have an irreducible character
$\chi\in \Irr(G)$ such that $\chi \overline{\chi}$ has at most two 
non-principal irreducible constituents.

\end{abstract}

\maketitle

\begin{section}{Introduction}

Let $G$ be a finite group. Denote by $\Irr(G)$ the set of 
irreducible complex characters of $G$. Let 
$1_G$ be the principal character of $G$. Denote by 
$[\Theta, \Phi]$ the inner product of the characters
$\Theta$ and $\Phi$ of $G$. Through this work, we will use the notation 
of \cite{isaacs}.

Let $\chi \in \Irr(G)$. Define $\overline{\chi}(g)$ to be
the complex conjugate $\overline{\chi(g)}$ of $\chi(g)$ 
for all $g \in G$. Then $\overline{\chi}$ is also 
an irreducible complex character of $G$. Since the product of 
characters is a character, then $\chi\overline{\chi}$ is a 
character of $G$. So it can be expressed as an integral linear
combination of irreducible characters. 
Now observe that 
	$$[\chi\overline{\chi}, 1_G]=[\chi, \chi]=1$$
\noindent where the last equality holds 
 since $\chi \in \Irr(G)$.
Assume now that $\chi(1)>1$. Then the decomposition of
the character $\chi \overline{\chi}$ into its irreducible 
constituents $\alpha_1, \ \alpha_2, \ldots, \alpha_n$
has the form
\begin{equation*}
 \chi\overline{\chi}= 1_G + \sum_{i=1}^n m_i \alpha_i
\end{equation*}
\noindent where
$m_i>0$ is the multiplicity of 
$\alpha_i$. Set $\eta(\chi)=n$. Thus $\eta(\chi)$ is the number of distinct non-principal
irreducible constituents of $\chi\overline{\chi}$.

In Theorem A of \cite{edithd}, it is proved that there exist universal constants
 $C$ and $D$ such that for any  finite solvable group $G$
 and any $\chi\in \Irr(G)$, we have that 
 $\dl(G/\Ker(\chi))\leq C\eta(\chi\overline{\chi})+D$. In Theorem B of 
 \cite{edithd}, it is proved that if $G$ is a  finite solvable group and
  $\chi\in \Irr(G)$, then $\chi(1)$ has at most $\eta(\chi)$ distinct
  prime divisors. If, in addition, $G$ is supersolvable, then 
  $\chi(1)$ is the product of at most $\eta(\chi)-1$ primes.
The purpose of this note is to study the solvable groups $G$ that have
faithful character $\chi\in \Irr(G)$ 
such that $\eta(\chi)\leq 2$.

Let $G$ be a solvable group and $\chi\in \Irr(G)$ be a faithful character.
 Assume that $\eta(\chi)=1$, i.e $\chi\overline{\chi}=1_G+m\alpha$ for 
 some $\alpha\in \Irr(G)$ and some integer $m>0$.
 Since $\chi\overline{\chi}$
and $1_G$ are real characters, so is $\alpha$. Thus $G$ has an irreducible real character.
Therefore $G$ has even order. More can be said in that regard and 
Theorem A ``classifies" those groups.

	\begin{teoremaa} 
Let $G$ be a finite solvable group and
 $\chi \in \Irr(G)$ be a faithful 
character. Assume that 

\begin{equation}\label{basicone}
\chi \overline{\chi} = 1_G + m\alpha
\end{equation}

\noindent for some $\alpha \in \Irr(G)$ and some 
positive integer $m$. Set $Z={\bf Z}(G)$. 
Let $E/Z$ be a chief factor of $G$. 
Set $\bar{G}=G/E$.
Then $m=1$, $\chi(1)^2 =|E:Z|$ 
 and one of the following holds:

(i) $\chi(1)=2$ and  $G/Z \cong A_4$ or $G/Z\cong S_4$. 

(ii) $\chi (1)=3 $ and either
$\bar{G}\cong  {\bf Q}_8$ or $\bar{G} \cong \SL(2,3)$.

(iii) $\chi(1)=5$
and $\bar{G}\cong \SL(2,3)$.
 
(iv) $\chi (1)=7$ and 
$\bar{G} \cong \widetilde{\GL}(2,3)$.

(v) $\chi(1)= 9$, ${\bf F}(\bar{G})$ is an
extra-special group of order $2^5$, namely the central product
of  the quaternion group
 ${\bf Q}_8$ and the dihedral group of order 8,
 $|{\bf F}_2(\bar{G})/{\bf F}(\bar{G})|=5$
and $\bar{G}/{\bf F}_2(\bar{G})$ is a 
subgroup of the cyclic group of order 2.
\end{teoremaa}

 We can check  that in all the cases given by Theorem A we have that
 $|G:Z|\leq 25920$. Assume now that $G$ is a solvable group, 
 $\chi\in \Irr(G)$ is a faithful character and   $\eta(\chi)=2$.
 In such situation, the index of the center $Z$ of $G$ in 
 $G$ is not longer bounded 
 by a universal constant, not even if in addition the group is nilpotent 
 (see Subsection 4.4).
 But the derived length $\dl(G)$ is at most 18 and we have 
 more information regarding the structure of the group $G$.

\begin{teoremab}
Let $G$ be a finite solvable group and 
$\chi \in \Irr(G)$ be a faithful character.
 Assume
 
\begin{equation*}
\chi\overline{\chi}= 1_G + m_1\alpha_1 + m_2\alpha_2
\end{equation*}
\noindent  where $\alpha_1, \alpha_2 \in \Irr(G)$ are non-principal
characters, $m_1$ and $m_2$ are strictly positive integers.
Let $Z$ be the center of $G$.
Then 
	
(i) The order of $G$ is even.

(ii) $\dl(G) \leq 18$.

(iii) Either $\Ker(\alpha_1)=Z$ or
$\Ker(\alpha_2)=Z$.
 If $\, \Ker(\alpha_1)$ and $\Ker (\alpha_2)$ are both abelian subgroups, then
$\chi(1)$ is a power of a prime number. Otherwise
$\{m_1,m_2\}=\{1,4\}$ and $\chi(1)\in \{10,14\}$. 	 
\end{teoremab}		

Propositions \ref{kernelsequal}, \ref{abelian} and 
 \ref{knotabelian} describe with more detail the possible structure of the group $G$.
 Examples of groups and characters satisfying the hypotheses of Theorem B are studied in 
 Subsection 4.3.
 
{\bf Acknowledgment.}
 This is part of my Ph.D. Thesis.
I thank Professor Everett C. Dade, my adviser, for his advise, suggestions and 
generosity with his time.
 I thank Professor I. Martin Isaacs for
 suggesting the question that started my studies in products of characters.
 I would also like to thank the mathematics department of the University 
of Wisconsin, at Madison, for 
their hospitality while I was visiting Professor Isaacs,
 and the mathematics department
of the University of  Illinois at Urbana-Champaign for their support. 

\end{section}

\begin{section}{Vector spaces}

B. Huppert classified
 all the solvable groups that act
faithfully on a finite vector space $V$ and transitively on 
$V^{\#}=V\setminus \{\,0\}$ (see Lemma \ref{huppert}).
 If in addition $V$ is a symplectic
vector space and  the action of $G$ on $V$ preserves the symplectic form,
 the structure of the group $G$ must be very 
constrained. We will classify the solvable
groups $G$ such that
$G$ can act faithfully on a symplectic vector 
space preserving the symplectic form 
and  acting transitively on $V^{\#}$.

In  Subsection  2.1 
we  mention the results that we  use
in the proof of Theorem \ref{symplectictrans}, the main result of this
section. 

\begin{subsection}{ Preliminaries }

	Let $V$ be a vector space of dimension $n$ over $\GF(q)$.
We may assume that $V$ is $\GF(q^n)$, the Galois field with
$q^n$ elements,  since $\GF(q^n)$ is  a vector
space of dimension $n$ over $\GF(q)$. Let 
$\mathcal{G}= \Gal(\GF(q^n)/ \GF(q))$ be the 
group of automorphism of the field $\GF(q^n)$ fixing
$\GF(q)$. Then  the \emph{semi-linear group} $\Gamma(V)$ of $V$ is 
	$$\Gamma(V)= \{ x \mapsto ax^{\sigma} \ | a \in (\GF(q^n))^{\#},
\sigma \in \mathcal{G}\}$$
Observe that 
\begin{equation}\label{unitsofgamma}
 U= \{ x \mapsto  ax \ | a \in (\GF(q^n))^{\#}\}
\end{equation}
\noindent  is a normal 
subgroup of $\Gamma(V)$. Also observe that $U \cong  (\GF(q^n))^{\#}$, the
group of units of the field $\GF(q^n)$.
So $U$ is a cyclic subgroup of $\Gamma(V)$ of order
$q^n-1$. Recall that
$\mathcal{G}$ is a cyclic group of order $n$.  
 Observe that $\Gamma(V)/U\cong \mathcal{G}$. Therefore
$\Gamma(V)/U$ is a cyclic group of order $n$.

	\begin{lemma}[B. Huppert]\label{huppert}
Let $V$ be a vector space of dimension $n$ over $\GF(q)$, where $q$
 is a prime power.
 Suppose that $G$ is a solvable subgroup of $GL(V)$
that transitively permutes the elements of $V^{\#}$. Let
$ {\bf F}(G)$ be the Fitting subgroup of $G$. 
Then one of the 
following holds:

       (i) $G$ is isomorphic to a subgroup of $\Gamma(V)$. 

	(ii) $q^n =3^4$, $ {\bf F}(G)$  is extra-special of order $2^5$,
$|{\bf F}_2(G)/ {\bf F}(G)|=5$ and $G/{\bf F}_2(G) \leq {\bf Z}_4$.

	(iii) $q^n= 3^2$, $5^2$, $7^2$, $11^2$ or 
$23^2$. Here ${\bf F}(G)= QT$, where $T= {\bf Z}(G)\leq {\bf Z}(GL(V))$
is cyclic, ${\bf Q}_8 \cong Q \trianglelefteq G$, $T\cap Q= {\bf Z}(Q)$
and $ Q/Z(Q)\cong {\bf F}(G)/T$ is a faithful irreducible
 $G/{\bf F}(G)$-module. We also have one of the following entries:

\begin{table}[ht]\label{table}
\begin{center}
\begin{tabular}{|c|c|c|} \hline
$q^n$ & $|T|$                 & $G/{\bf F}(G)$  \\ \hline
$3^2$ &  2                    & ${\bf Z}_3 {\text \ or\ } S_3$ \\ \hline
$5^2$ & 2 or 4           & ${\bf Z}_3$ \\ \hline
$5^2$ & $4$                     &  $S_3 $\\  \hline
$7^2$ & 2  or 6                & $S_3 $\\  \hline
$11^2$ & $10$                   & ${\bf Z}_3 {\text \ or\ } S_3 $ \\ \hline
$23^2$ & $22$                  & $S_3$ \\ \hline
\end{tabular}
\end{center}
\end{table}

\end{lemma}
\begin{proof}
See Theorem 6.8 of \cite{wolf}.
\end{proof}


	Let $V$ be a finite dimensional vector space over the field 
$F$.
We say that $V$ is a \emph{symplectic vector space} if there exists a 
non-degenerate, alternating and bilinear form $f$ on $V$.
	
	We denote by $\Spl(V)$ the group of automorphism of $V$
preserving the symplectic form $f$. We denote by 
$\Spl(n, q)$ the group of automorphism of a $n$-dimensional
symplectic vector space over $\GF(q)$.
The following is a list of well known results that we will use later.
 
\begin{lemma}\label{symplectic2}
 Let $\Spl(2,q)$ be the symplectic group where
$(q,2)=1$.
Then 

 (i) $\Spl(2,q) =  \SL(2,q)$.

(ii) $|{\bf Z} (\SL (2,q))| = 2$.

(iii) If $q> 3$, then $\SL(2, q)$ is not a solvable group.

 (iv) $| \SL (2,q)| = q(q^2-1)$.

(v) $\SL(2,3)\cong {\bf Q}_8 \rtimes {\bf Z}_3$.

 (vi) $\SL(2,3)$ is isomorphic to a subgroup of $\SL(2,5)$.

(vii) The 2-Sylow subgroups of $\GL(2,3)$ are isomorphic to the semi-dihedral
group of order $16$.

(viii) For $q$ odd, 
the 2-Sylow subgroups of $\SL(2,q)$ are isomorphic to 
generalized quaternion groups.
\end{lemma}

\begin{definition}\label{twin}
 We define $\widetilde{\GL}(2,3)$ as the group 
given by the following representation
\begin{equation*}
< x, y , z \mid  x^8=1, \  y^2= x^4, \ z^3=1, \  x^y =x^{-1},
 \ z^x= x^2 y z^{-1}, \   (x^2)^z=y>.
\end{equation*}
\end{definition}

\begin{lemma}
$\widetilde{\GL}(2,3)$ is not isomorphic to $\GL(2,3)$.
Also $\widetilde{\GL}(2,3)$ has a subgroup of index 2 isomorphic 
to $\SL(2,3)$, and ${\bf F}(\widetilde{\GL}(2,3))\cong {\bf Q}_8$.
\end{lemma}
\begin{proof}
We can check that the subgroup generated by 
$x$ and $y$ is a 2-Sylow subgroup of $\widetilde{\GL}(2,3)$ and 
is  
isomorphic to ${\bf Q}_{16}$. Thus $\GL(2,3)$ and 
 $\widetilde{\GL}(2,3)$ are not isomorphic since 
a 2-Sylow subgroup of $\GL(2,3)$ is isomorphic to the 
semi-dihedral group of order 16.

Observe that 
${\bf F}(\widetilde{\GL}(2,3)) ={\bf Q}_8$, namely the subgroup 
generated by $x^2$ and $y$. Also observe that 
 the subgroup 
generated by $x^2$, $y$ and $z$ is isomorphic to $\SL(2,3)$.
\end{proof}
The group $\widetilde{\GL}(2,3)$ is the isoclinic group of 
$\GL(2,3)$. 
\begin{lemma} \label{sl27}
The group $\SL(2,7)$ has a subgroup isomorphic to $\widetilde{\GL}(2,3)$.
Also, any subgroup of order 48 in $\SL(2,7)$ is isomorphic to 
$\widetilde{\GL}(2,3)$.
\end{lemma}
\begin{proof}

Let $x, y, z \in \SL(2,7)$ be the matrices
\begin{equation*}
x=\left( \begin{array}{cc} 2 & 2 \\ 5 & 2 \end{array} \right), \ \ 
 y=\left( \begin{array}{cc} 3 & 2 \\ 2 & 4 \end{array} \right), \ \ 
 z=\left( \begin{array}{cc} 4 & 0 \\ 1 & 2 \end{array} \right).
\end{equation*} 
We can check that they satisfy the relations given 
in Definition \ref{twin}. Thus $\SL(2,7)$
has a subgroup isomorphic to $\widetilde{\GL}(2,3)$, mainly
the group generated by $x$, $y$ and $z$.  

Let $H$ be a subgroup of $\SL(2,7)$ of order 48. Then $H$ is a maximal
subgroup of $\SL(2,7)$.  We can check 
using the Atlas  that  the maximal subgroups of $\SL(2,7)$ of
order 48 are all isomorphic.
\end{proof}

\end{subsection}
\begin{subsection}{Transitive actions and symplectic vector spaces}

\begin{theorem}\label{symplectictrans}
Let $G$ be a finite solvable group. Let $V$ be a  
symplectic vector space of dimension $n$ over $\GF(q)$, where $q$ is
a power of a  prime number.
 Assume that $V$ is a faithful
$G$-module
and the action of $G$ on $V$ preserves the symplectic form.
Also assume that  $G$ acts transitively on $V^{\#}$.
 Set $e^2=|V|=q^n$.
Let ${\bf F}(G)$ be the Fitting subgroup of $G$ and
 ${\bf F}_2(G)/{\bf F}(G)$
be the Fitting subgroup of $G/{\bf F}(G)$.
Then one of the following holds:
	
	(i) $e=2$ and   $G\cong {\bf Z}_3$ or $G\cong S_3$.

	(ii) $e=3$, and 
either $G\cong {\bf Q}_8$ or $G\cong \SL(2,3)$.

         (iii) $e=5$ and
$G\cong \SL(2,3)$.

	(iv) $e=7$ and 
$G\cong \widetilde{\GL}(2,3)$.

	(v) $e=9$, ${\bf F}(G)$ is an
extra-special group of order $2^5$, namely the central product
of  the quaternion group
 ${\bf Q}_8$ and the dihedral group of order 8, 
$|{\bf F}_2(G)/{\bf F}(G)|=5$
and $G/{\bf F}_2(G)$ is a 
subgroup of the cyclic group of order 2.
\end{theorem}

Observe that the hypothesis of Theorem \ref{symplectictrans}
 implies the hypothesis
of Lemma \ref{huppert}. Thus either (i), (ii) or (iii) of 
Lemma \ref{huppert} holds.
 We will analyze each case given by the conclusion of 
Lemma \ref{huppert}.

\begin{lemma} \label{symplectic}

Let $C$ be a cyclic group, $|C|=n$. Let $V$ be a finite
symplectic vector space.
Assume $C$ acts on $V$ preserving the symplectic form.
Assume that $V$ is an irreducible and faithful $C$-module.
Set $|V|=e^2$.
Then $n$ divides $e+1$.
\end{lemma}

\begin{proof} See \cite{huppert}, Satz 9.23.
\end{proof}
We say that a group $G$ acts \emph{Frobeniusly} on a vector space $V$ if 
every nontrivial element of $G$ acts fixed-point-free on $V$.

\begin{lemma}\label{cyclic} Let $V$ be a finite 
symplectic vector space.
Let $G$ be a finite solvable group that acts faithfully, irreducibly
 and symplectically 
on $V$. Set $|V|=e^2$.
 Assume ${\bf F}(G)$ acts Frobeniusly on $V$.
 Let $M$ be a normal cyclic subgroup of ${\bf F}(G)$. Then $|M| \leq e+1$. 
\end{lemma}

\begin{proof}
Set $F={\bf F}(G)$.
Since $F$ is normal in $G$, by Clifford theory 
we have that $V$ decomposes as a direct sum of
irreducible $F$-modules, each of the same size.
	
	Our hypothesis is that $F$ acts Frobeniusly on $V$. Assume that
$V$ is not an irreducible $F$-module. By Clifford theory 
we have that $V$ decomposes as a direct sum of at least 2
irreducible $F$-modules, each of the same size. Since $F$ acts
Frobeniusly on $V$, we have that $|F|\leq |U|$, where 
$U$ is one of the irreducible components of the $F$-module $V$.
Thus $|F|\leq e$. Since $M$ is a subgroup of $F$, it follows
that $|M|\leq e< e+1$. 

We may assume that $V$ is an irreducible $F$-module.
By hypothesis, we have that $V$ is a symplectic space and that $G$ acts on 
$V$ preserving the symplectic form. If $V$ is not an irreducible
$M$-module, then as before we can conclude that $|M|\leq e$.
So we may assume that $V$ is an irreducible $M$-module.
Then by Lemma \ref{symplectic}
we have that $|M|$ divides $e+1$. 
\end{proof}

\begin{lemma}\label{sizef}
Assume $G$ is a finite solvable group.
Then 
\begin{equation*}
|G| \mbox{ divides }|\Aut({\bf F}(G))||Z({\bf F}(G))|.
\end{equation*}
\end{lemma}

\begin{proof}[Proof of Theorem \ref{symplectictrans}]
	Since the  hypothesis of Theorem \ref{symplectictrans} 
implies the
hypothesis of Lemma \ref{huppert}, we will study the cases
given by the conclusion of Lemma \ref{huppert}.

\begin{claim}\label{tied} $e^2-1=q^n-1$ divides the order of $G$.
\end{claim}

\begin{proof} By Hypothesis, $G$ acts
transitively on $V^{\#}$ and $e^2= |V|$.
 Therefore $e^2-1$ divides the order of $G$.
\end{proof}

\begin{claim}\label{g/hactsh}
 Assume that $G$ is a subgroup of the semi-linear
group $\Gamma(V)$. Also assume that $G$ acts transitively
on $V^{\#}$. Set $H= G \cap U$, where $U$ is as in 
\eqref{unitsofgamma}. Then $H$ is a normal cyclic subgroup of
$G$ and $G/H$ acts faithfully on
$H$. 
\end{claim}
\begin{proof} 
Since $U$ is a normal cyclic subgroup of $\Gamma(V)$, 
we have that
$H= U\cap G$ is a normal cyclic subgroup of $G$. Since
$\Gamma(V)/U$ is a cyclic group of order $n$, $G/H$ is a 
cyclic group of order a divisor of $n$.

Since $H$ is a cyclic normal subgroup of $G$, we have that
$G/H$ acts on 
$H$ by conjugation. 
Suppose that $G/H$ does not act faithfully on $H$.   
Let $K/H$ be the kernel of the action and $t= |K/H|$.
Thus $G/K$ acts faithfully on $H$. Since $H$ is cyclic,
$|G/K|$ divides $\varphi(|H|)$, where $\varphi$ is 
the Euler function. Since $K/H$ is the kernel of the action of
$G/H$ on $H$, $t= |K/H|$ and $G/H$ is isomorphic
to a subgroup of $\mathcal{G}$, we have that $H$ is isomorphic
to a subgroup of  the group of units $(\GF(q^{n/t}))^{\#}$ 
of a subfield
of $\GF(q^n)$. Thus $|H|$ divides $q^{n/t}-1$. Therefore we 
have that
\begin{equation}\label{orderg}
|G| \mbox{ divides } (q^{n/t}-1)t \varphi(q^{n/t}-1).
\end{equation}
Observe that
 $q^n-1= (q^{n/t}-1)(q^{\frac{(t-1)n}{t}}+ q^{\frac{(t-2)n}{t}}+\cdots + 1)$.
Thus by \eqref{orderg} and Claim \ref{tied} we have that
\begin{equation}\label{goingcarefully}
(q^{\frac{(t-1)n}{t}}+ q^{\frac{(t-2)n}{t}}+\cdots + 1) \mbox{ divides }
t \varphi(q^{n/t}-1).
\end{equation}
Observe that $1+q^{\frac{2n}{t}}>2q^{n/t}$ since $(1-q^{n/t})^2>0$.
If $t>2$, we have that 
$$q^{\frac{(t-1)n}{t}}+ q^{\frac{(t-2)n}{t}}+\cdots + 1> tq^{n/t}> 
t \varphi(q^{n/t}-1).$$
\noindent That can't hold by \eqref{goingcarefully}. Thus $t=2$.
If $q$ is odd, then $q^{n/2} -1$ is even. That implies 
$$\varphi(q^{n/2}-1)\leq \frac{1}{2} (q^{n/2}-1).$$
\noindent But then $$q^{n/2}+1> q^{n/2}-1 \geq 2\varphi(q^{n/2}-1).$$
By \eqref{goingcarefully} we conclude that $q$ must be even.
Since $\varphi (q^{n/2}-1) \leq q^{n/2} -1$, 
 by \eqref{goingcarefully} we have that 
 $$q^{n/2}+1\, \bigm| \, 2 \varphi(q^{n/2}-1) \leq 2(q^{n/2}-1).$$
Thus $q^{n/2}+1 =2\varphi(q^{n/2}-1)$ and $q$ can not be even. 
Therefore $t=1$ and so the action of $G/H$ on $H$ is faithful.
\end{proof}

\begin{case} Assume (i) of Lemma \ref{huppert}, i.e.,
assume $G$ is a subgroup of the semi-linear group $\Gamma(V)$. 
Then 

(i) $e=2$ and $G \cong {\bf Z}_3$ or  $G \cong  S_3$.

(ii) $e=3$ and $G\cong {\bf Q}_8$.
\end{case}
 
\begin{proof}
By Claim \ref{g/hactsh}, we have that $H= G \cap U$ is a normal
cyclic subgroup of $G$ and $G/H$ acts faithfully on $H$.
  
 Because $H$ is cyclic, we have 
then $|\Aut (H)| = \varphi(|H|)$, where $\varphi$ is the Euler
function. 
Since $G/H$ acts faithfully on $H$, we have that
\begin{equation*}\label{eq1}
|G|\leq |H| \varphi(|H|).
\end{equation*}

	If $|H|$ is not a prime number 
and $|H|\neq 4$ , then $\varphi(|H|) \leq |H|-3$.
By Lemma   \ref{cyclic}
  we have that $|H| \leq e+1$. So 
$$ \varphi(|H|) \leq |H|-3\leq
e+1-3 < e-2.$$
\noindent But then $|G|\leq |H|\varphi(|H|)< e^2-1$, a contradiction
with Claim \ref{tied}.
	So we can assume that either $|H|$ is a prime number or $|H|=4$.

	Suppose that $|H|$ is a prime number $q$. 
Then by Lemma \ref{cyclic}
we have that $q\leq e+1$. If $q< e+1$, then $q\leq e$ and 

\begin{equation*}
|G|\leq |H|\varphi(|H|)= q(q-1)\leq e (e-1) < e^2-1.
\end{equation*}
But that is a contradiction with Claim \ref{tied}. Thus $q=e+1$.
Observe that $q(q-2)= e^2-1$. Because of Claim \ref{tied},
we have that $q(q-2)= e^2-1$ divides $q(q-1)$ . Since $e>1$,
$q=3$ and 
$e=2$. Thus $|H|=3$ if $|H|$ is a prime number.
 
	It remains to consider $|H|=4$.
By Lemma \ref{sizef}
$|G|$ divides $|H|\varphi(|H|)= 8 $. By Claim  \ref{tied}, we
 have that $e^2-1$ divides $|G|$. Since $|G|$ divides 8,
then $e^2-1=|V^{\#}|$ is 2, 4 or 8. Since $V$ is a symplectic 
vector space,  $|V|=q^2$ for some power of a prime $q$. 
It follows that $|V^{\#}|=8$ and $G\leq \SL(2,3)$. Thus 
$G\cong {\bf Q}_8$.
\end{proof}	

\begin{case} Assume case (ii) of Lemma \ref{huppert}. Then

(v) $q^n =3^4$, $ {\bf F}(G)$  is extra-special of order $2^5$,
$|{\bf F}_2(G)/ {\bf F}(G)|=5$ and $G/{\bf F}_2(G) \leq {\bf Z}_2$.
\end{case}

\begin{proof} Since $G$ acts symplectically and faithfully on
a vector space $V$ and $|V|= 3^4$, we have that $G$ is isomorphic
to a subgroup of $\Spl(3,4)$.
Using the Atlas, we can check that  $G/{\bf F}_2(G) \leq {\bf Z}_2$.
\end{proof}
\begin{case} Assume case (iii) of Lemma \ref{huppert}.
Then $e \in \{5,7\}$. Furthermore one of the following holds
	
	(ii) $e=3$ and $G \cong  \SL(2,3)$.

        (iii) $e=5$ and $G \cong \SL(2,3)$.

	(iv) $e=7$ and $G \cong \widetilde{\GL}(2,3)$.
\end{case}

\begin{proof} 
By Lemma \ref{huppert} (iii), we have that 
$e^2\in \{3^2, 5^2, 7^2, 11^2, 23^2\}$. Thus $G \leq \Spl(2,e)$. 
Recall that ${\bf F}(G)= {\bf Q}_8 T$, where  $T \leq {\bf Z}(\GL(V))$. 
Thus $T=  {\bf Z}(\GL(V)) \cap \Spl(2,e)= {\bf Z}(\Spl(2,e))$.
Since $\Spl(2,e)= \SL(2,e)$ and 
$|{\bf Z}(\SL(2,e))|=2$ by 
 Lemma \ref{symplectic2} (i) and (iii),
we have that $|T|=1$ or $|T|=2$. Since $T \cap {\bf Q}_8= {\bf Z}({\bf Q}_8)$,
it follows that $T= {\bf Z}({\bf Q}_8)$.
Thus ${\bf F}(G)={\bf Q}_8$.

Since $\Aut({\bf Q}_8) \cong S_4$ (see Exercise 5.3.4 of \cite{robinson}),
we have that   $|\Aut({\bf Q}_8)|=24$. Since ${\bf Z}({\bf Q}_8)$
 is cyclic of order 2, 
by Claim \ref{sizef} we have that $|G|$ is a 
divisor of $48$. By Claim \ref{tied}
we have that $e^2-1$ has to divide 48.
Therefore $e\in \{3,5,7\}$.

	By hypothesis $G$ is isomorphic to a subgroup of
$\Spl(V)$. Recall that  $|V|=e^2$. By Lemma \ref{symplectic2}
we have that $\SL(2,e)=\Spl(2,e)$. Thus if
$e=3$ then $|G|$ divides $|\SL(2,3)|=24$. Since $G$ acts
transitively on a set of $8$ elements, we have that
8 divides $|G|$.  Since 
 $ {\bf Q}_8/Z({\bf Q}_8)\cong {\bf F}(G)/T$ is a faithful irreducible
 $G/{\bf F}(G)$-module, it follows that $G \cong \SL(2,3)$.

If $e=5$, then $G$ is isomorphic to a solvable subgroup  	 
of $\SL(2,5)$. By Claim \ref{tied} we have that 
 the order of $G$ is divisible by $24=e^2-1$.
Since  
$\SL(2,5)$ is not solvable (see Lemma \ref{symplectic2} (iii) ) 
and $|\SL(2,5)|=5\times 24$,
we have that $|G|= 24$. Since ${\bf F}(G) \cong {\bf Q}_8$,
and  $|G|= 24$, we have that $G \cong \SL(2,3)$.

Similarly, if $e=7$, $G$ is isomorphic to a subgroup of 
$\SL(2, 7)$. By Claim \ref{tied} we have that
 $G$ is a solvable group of order
divisible by $48=7^2-1$.
Since $\SL(2,7)$ is not a 
solvable group (see Lemma \ref{symplectic2} (iii) ) and
$|\SL(2,7)|= 7\times 48$ we have that 
$|G|=48$. By Lemma \ref{sl27}, we have that 
$G\cong \widetilde{\GL}(2,3)$.
\end{proof}

We have analyzed all the cases given 
in the conclusion of Lemma \ref{huppert}. 
So we have finished the proof.
\end{proof}

\end{subsection}

\end{section}

\begin{section}{Theorem A}

Our main result of this section is Theorem A.
In Proposition \ref{allexist}, for each case (i) to (v)
 in the conclusion
of Theorem A, we will give an 
example of a solvable group and a faithful character 
$\chi \in \Irr(G)$ satisfying the conditions in that case.

\begin{subsection}{General results}

All the  results in this section are proved in \cite{edithd}.

\begin{lemma} \label{basico1} Assume that  $G$ is a finite group.
Let $L$ and $N$ be normal subgroups of $G$ such that
$L/N$ is an abelian chief factor of $G$. Let $\theta \in \Irr(L)$
be a $G$-invariant character such that $\theta_N$ is reducible. 
	Then 
\begin{equation*}\label{principalcharacter}
\theta\overline{\theta}= 1_N^L + \Phi
\end{equation*}
\noindent where $\Phi$ is either the zero function or a character of $L$ and $[\Phi_N, 1_N]=0$.

\end{lemma}
\begin{proof} See Lemma 4.1 of \cite{edithd}.
\end{proof}
\begin{lemma}\label{center}
Assume that $G$ is a finite group and $\chi \in \Irr(G)$ is a faithful character. 
Let
$\{\alpha_i \in \Irr(G)^{\#}\mid i=1, \ldots , n \}$
be the set of 
non-principal irreducible constituents of $\chi \overline{\chi}$.
 Then
\begin{equation*}\label{centerequ}
{\bf Z}(G)= \bigcap_{i=1}^n \Ker(\alpha_i).
\end{equation*}
\end{lemma}
\begin{proof}
See Lemma 5.1 of \cite{edithd}.
\end{proof}
\begin{lemma}\label{a1equal1} 
Assume that  $G$ is a finite solvable group and $\chi \in \Irr(G)$ 
with $\chi(1)>1$.
Let $\{\alpha_i \in \Irr(G)^{\#}\mid i=1, \ldots , n \}$
be the set of 
non-principal irreducible constituents of $\chi \overline{\chi}$.
If   
 $\, \Ker(\alpha_j)$
 is not properly contained in $\Ker(\alpha_i)$
for all $i$, then $ [ \chi \overline{\chi} , \alpha_j] =1$. 
Thus  $1\in \{ [\chi \overline{\chi} , \alpha_i] \mid 
 i=1, \ldots, n\}$.

\end{lemma}
\begin{proof} See Theorem C of \cite{edithd}.
\end{proof}
\end{subsection}

\begin{subsection}{Proof of Theorem A}

We will need two lemmas in the proof of 
Theorem A. These results are well known.

\begin{lemma}\label{centralizere}
 Let $G$ be a solvable group with cyclic center ${\bf Z}(G)$. 
 Assume that
$E/{\bf Z}(G)$ is a chief factor of $G$. If $E/{\bf Z}(G)$
is a fully ramified
section with respect to some $\Lambda \in \Irr(E)$, 
then ${\bf C}_G(E/{\bf Z}(G))=E{\bf C}_G(E)$. 
Also $E/{\bf Z}(G)$
is a symplectic vector space and the conjugate
action of $G$ on $E/{\bf Z}(G)$ preserves the symplectic form.
\end{lemma}

\begin{lemma}\label{restrictionirreducible}
Let $G$ be a group and $H$ be a subgroup of $G$. If $\chi \in\Irr(G)$ 
and $\chi_H \in \Irr(H)$, then ${\bf C}_G(H) \leq {\bf Z}(\chi)$.
\end{lemma}

\begin{proof}[Proof of Theorem A]

Lemma \ref{a1equal1} implies that $m=1$.

Since $\chi$ is a faithful character, we have that 
 ${\bf Z}(\chi)={\bf Z}(G)=Z$.
By Lemma \ref{center} we have that  $ \Ker(\alpha)=Z$.

\begin{claim}\label{normal}
 Let $N$ be a normal subgroup of $G$. Then $\chi_N \in \Irr(N)$ 
or $N \leq Z$.
\end{claim}

\begin{proof}
Assume $N$ is not contained in $Z$. Since 
$\Ker(\alpha)= Z$ and $N$ is normal in $G$, 
we have that $[\alpha_N, 1_N]=0$.  
We conclude that $[(\chi {\overline{\chi}})_N, 1_N]=1$.
Thus $\chi_N$ is an irreducible character of $N$.
\end{proof}

\begin{claim}\label{fully} 
Let $E/Z$ be a chief factor of $G$. Let $\lambda \in \Irr(Z)$ be
such that $[\lambda, \chi_Z]\neq 0$. Then $\chi_E$ and $\lambda$
are fully ramified with respect to $E/Z$. Therefore $E/Z$ is a symplectic
vector space and $G$ acts on it preserving the symplectic form.
Also $\chi(1)^2=|E:Z|$.
\end{claim}

\begin{proof}
By Claim \ref{normal}, we have that  $\chi_E \in \Irr(E)$.
By \eqref{basicone} we have that $$\chi(1)^2= 1 + m \lambda(1)$$ 
\noindent   
 and therefore $\chi(1)>1$.
 By hypothesis, $G$ is solvable
and so $E/Z$ is abelian. Since $Z={\bf Z}(G)$ and $\chi_Z$ reduces,
by Theorem 6.18 of \cite{isaacs} we have that $\chi_E$ is fully ramified
with respect to $E/Z$. Thus
$\chi_Z = e \lambda$, where $ e^2 =|E:Z|$.
In particular, we have that $e= \chi(1)$.
\end{proof}

\begin{claim}\label{a1}
 $(\chi\overline{\chi})_E= 1_Z^E$.
\end{claim}

\begin{proof} Since $\chi(1)^2=|E:Z|$, by Lemma \ref{basico1} we get 
the claim. 
\end{proof}

\begin{claim}\label{transitive} 
$C_G(E/Z)=E$. Thus $G/E$ acts symplectically and faithfully
 on $E/Z$, transitively on
${(E/Z)}^{\#}$.
\end{claim}
\begin{proof}
By Lemma \ref{centralizere} we have that 
${\bf C}_G(E/Z)= E {\bf C}_G(E)$. If ${\bf C}_G(E) \neq Z$, then 
Claim \ref{normal} implies that 
$\chi_{{\bf C}_G(E)} \in \Irr({\bf C}_G(E))$.
By Lemma \ref{restrictionirreducible} we have that 
\begin{equation*}
E\leq {\bf C}_G({\bf C}_G(E))\leq Z(\chi) ={\bf Z}(G)=Z<E.
\end{equation*}
We conclude that
${\bf C}_G(E)= Z$. Thus $E = {\bf C}_G(E/Z)$.
Therefore $G/E$ acts faithfully on $E/Z$.

By Claim \ref{a1} we have that $\alpha_E=1_Z^E  - 1_E$. 
Since $E$ is normal in $G$ and $\alpha \in \Irr(G)$, by 
Clifford theory  we have that
$G$ acts transitively on the set 
$$\Irr(E\modu Z)^{\#}=
\{ \gamma \in \Irr(E) \mid  \gamma \neq 1, \Ker(\gamma)\geq Z\}.$$
\noindent Therefore by  Lemma 5.2 of \cite{go}, 
we have that $G$ acts 
transitively
on  ${(E/Z)}^{\#}$. The rest follows from Claim \ref{fully} 
\end{proof}
\begin{claim}\label{symplecticoptions}
 Set $\bar{G}=G/E$. Then one of the following holds

	(i) $e=2$ and either
  $G\cong S_3$ or $G\cong {\bf Z}_3$.

	(ii) $e=3$ and 
either $\bar{G}\cong {\bf Q}_8$ or $\bar{G}\cong \SL(2,3)$.

         (iii) $e=5$ and
$\bar{G}\cong \SL(2,3)$.

(iv) $e=7$  and 
$\bar{G}\cong \widetilde{\GL}(2,3)$.

	(v) $e=9$, ${\bf F}(\bar{G})$ is an
extra-special group of order $2^5$,
namely the central product
of  the quaternion group
 ${\bf Q}_8$ and the dihedral group of order 8,
$|{\bf F}_2(\bar{G})/{\bf F}(\bar{G})|=5$
and $G/{\bf F}_2(\bar{G})$ is a 
subgroup of the cyclic group of order 2.
\end{claim}
\begin{proof}
By Claim \ref{transitive} we have the hypotheses of 
 Theorem \ref{symplectictrans}. The rest follows from 
Theorem \ref{symplectictrans}.
\end{proof}
Since $e^2= |E/Z|$ and $\chi(1)=e$, Theorem A follows from
Claim \ref{symplecticoptions}.
\end{proof}
\end{subsection}
\begin{subsection}{Examples}

The following is a set of examples of solvable groups with 
irreducible characters that satisfy the hypotheses of 
Theorem  A. We will use the notation of
Theorem A in this subsection.
 We will start with some lemmas that we will use
 in the construction
of some of the examples.  

\begin{lemma}\label{splitting}
 Let $E$ be an extra-special group of order $p^{2n+1}$ and 
exponent $p$,  for some odd prime $p$.
Observe that $E/{\bf Z}(E)$ is a vector space of dimension 
$2n$ over $\GF(p)$ with a
symplectic form. Then the following short exact sequence
\begin{equation*}
0 \rightarrow \Inn (E) \rightarrow \Aut(E\modu {\bf Z}(E))\rightarrow
 \Spl(2n,p)\rightarrow 0
\end{equation*}
\noindent  splits.
\end{lemma}

\begin{proof}
See \cite{winter}
\end{proof}

\begin{lemma}\label{itworksifchiextension} 
Let $G$ be a finite group. Let 
$E/{\bf Z}(G)$ be  an abelian 
 chief factor of $G$. 
Assume that the action by conjugation of 
$G/E$ on $(E/{\bf Z}(G))^{\#}$ is transitive. 
Let $\chi \in \Irr(G)$ be a faithful 
character of
$G$. Assume that 
$\chi_E \in \Irr(E)$ and $\chi(1)>1$.
Then $\chi$ satisfies \eqref{basicone},
 i.e.
$\eta(\chi)=1$.
\end{lemma}

\begin{proof}
Set $Z={\bf Z}(G)$.
Let $\lambda \in \Irr(Z)$ be a character of $Z$
such that $[\chi_Z, \lambda]\neq 0$. Since $\chi_E\in \Irr(E)$
and $Z={\bf Z}(G)\leq {\bf Z}(E)$, by Corollary 2.30 of 
\cite{isaacs} we have that
$\chi(1)^2\leq |E:Z|$. 
Observe that  
$\chi_Z$ is a reducible character since $\chi(1)>1$ and $Z={\bf Z}(G)$ is
abelian. Observe also
that $ 1_Z^E(1) = |E:Z| \geq \chi(1)^2$.
Since $\chi(1)^2\leq |E:Z|$,
by Lemma \ref{basico1} we conclude that
\begin{equation}\label{isfullyramified}
(\chi \overline{\chi})_E = 	 1_Z^E.
\end{equation}
\noindent In particular, $\chi(1)^2=|E:Z|$.

Let $\alpha \in \Irr(G)^{\#}$ be such that
$[\alpha, \chi\overline{\chi}]\neq 0$.
Observe such 
$\alpha$ exists since $\chi(1)>1$. 
By Lemma \ref{center},
 we have that $\Ker(\alpha)\geq Z$. By \eqref{isfullyramified}
we have that $E$ is not contained in $\Ker(\alpha)$. 
Thus the irreducible constituents of $\alpha_E$ are elements
of $\Irr(E\modu Z)^{\#}$. By Clifford theory 
we have that the irreducible constituents of $\alpha_E$ are 
$G$-conjugate. Our  hypothesis is  that $G/E$ acts transitively on 
$(E/Z)^{\#}$. Thus $G/E$ acts transitively on the set 
$\Irr(E\modu Z)^{\#}$.
We conclude that 
the set of irreducible constituents of 
$\alpha_E $ is $\Irr(E \modu Z)^{\#}$.
Therefore $\alpha(1)\geq |E:Z|-1$.      
Since $\chi(1)^2=|E:Z|$, 
$\alpha$ is an  irreducible constituent  
of $\chi\overline{\chi}$, $\alpha(1)\geq |E:Z|-1$ and 
$[1_G, \chi\overline{\chi}]=1$,
we have that 
$$\chi\overline{\chi}= 1_G + \alpha.$$
\end{proof}
\begin{lemma}\label{caseeasier}
Assume the notation of Theorem A.
 Let $p\in \{3,5,7\}$.
 There exists a group $G$ that satisfies the hypotheses
of Theorem A and such that $\chi(1)=p$ and
 $|G:{\bf Z}(G)|=p^2\times  (p^2-1)$. In particular

(i) If $p=3$, then $\bar{G} \cong {\bf Q}_8$.

(ii) If $p=5$, then $\bar{G} \cong \SL(2,3)$.

(iii) If $p=7$, then $\bar{G} \cong \widetilde{\GL}(2,3)$. 

\end{lemma}

\begin{proof}
Let $E$ be an extra-special group of order $p^3$ and 
exponent $p$, where
$p \in \{3, 5,7\}$. 

Observe that $\SL(2,3)$ has a subgroup of order 8, namely the
quaternion group ${\bf Q}_8$. We can check that $\SL(2,5)$ 
has a subgroup
isomorphic to $\SL(2,3)$ (see Lemma \ref{symplectic2} (v)), 
and thus  has a subgroup of order 
24.   Also $\SL(2,7)$ has a subgroup of order 48, namely 
$\widetilde{\GL}(2,3)$ (see Lemma \ref{sl27}).
Thus for $p \in \{3,5,7\}$, there exists a solvable subgroup
$H$ of $\SL(2,p)$ of order $p^2-1$.

Set  $V=E/{\bf Z}(E)$.
 Observe that $V$ is a vector space of dimension 2 over 
$GF(p)$ and $V$ has a symplectic form, namely 
$(\overline{v}, \overline{u})\mapsto [v,u]$.
Denote by $\Spl(2,p)$ the group of automorphisms
of $V$ that preserves the symplectic form of $V$.
 Since $\SL(2,p)=\Spl(2,p)$ (see Lemma \ref{symplectic2} (i)),
we may assume that $H$ is a subgroup of $\Spl(2,p)$.
 By Lemma \ref{splitting} we may assume that 
$H\leq  \Aut(E \modu {\bf Z}(E))$.
We can check that 
$H$ acts transitively on $V^{\#}$.

	Let $G=H \ltimes E$ be the semidirect product of $E$ by $H$.
 Observe that ${\bf Z}(G)= {\bf Z}(E)$.
Set $Z={\bf Z}(G)$.
 Let $\Lambda \in \Irr(E)$ be a character of degree $p$.
Observe that $\Lambda$ is a $G$-invariant character. Since
$(|H|, |E|)=1$, and $\Lambda$ is a $G$-invariant character, 
 $\Lambda$ extends to $G$ (see Theorem 13.3 of \cite{isaacs}).
  Let $\chi \in \Irr(G)$ be an extension
of $\Lambda$. 

	Observe that $\Ker(\chi)$ is a normal subgroup 
of $G$ and acts on 
$E$ by conjugation. Thus $\Ker(\chi)$ is normal in $\Ker(\chi)E$.
Since $E$ is also normal in $\Ker(\chi)E$, we have  
$\Ker(\chi)\leq {\bf C}_G(E)\leq E $.   Since
$\Ker(\chi) \cap E=\Ker(\Lambda)=1$, it follows that $\chi$ is a 
faithful character. 

	Observe that we have the hypothesis of Lemma 
\ref{itworksifchiextension}. Thus the result follows from 
that Lemma.
\end{proof}
\begin{proposition}\label{allexist}
For each of the cases (i) - (v) in the conclusion of 
 Theorem A, there 
exists a solvable group and a faithful character $\chi \in \Irr(G)$
satisfying the conditions in that case.
\end{proposition}
\begin{proof}

{\bf (i)} 
 Let $G=\GL(2,3)$ and $\chi\in \Irr(G)$ be a character of degree 2. We can
check that $\chi$ satisfies \eqref{basicone}. We can check that
$G/ {\bf Z}(G) \cong S_4$.

Let $G=\SL(2,3)$ and $\chi\in \Irr(G)$ be a character of degree 2. We can
check that $\chi$ satisfies \eqref{basicone}. We can check that
$G/ {\bf Z}(G) \cong A_4$.

{\bf (ii)} 
Lemma \ref{caseeasier} gives us  a solvable group satisfying
 (ii) of Theorem A with $\bar{G} \cong {\bf Q}_8$. 
It remains to prove that there exist a group $G$ and a faithful
character $\chi$ satisfying  (ii) such that
$\bar{G} \cong \SL(2,3)$.

Let $E$ be an extraspecial group of order $27$ and exponent 3. Then 
$E/ {\bf Z}(E)=V$ may be regard as a vector space of dimension 
2 over  $\GF(3)$, the finite field with  3 elements. 
Observe that
$\SL(2,3)$ acts on $V$ transitively. 
By Lemma \ref{symplectic2} (i)
we have that $\SL(2,3)=\Spl(2,3)$. 
Also by Lemma \ref{symplectic2}
(v), we have that $\SL(2,3)\cong {\bf Q}_8 \rtimes {\bf Z}_3$. We will assume
that $\SL(2,3)= {\bf Q}_8 \rtimes {\bf Z}_3$  By
Lemma \ref{splitting} we may regard $\SL(2,3)$ as a subgroup
of $\Aut(E \modu {\bf Z}(E))$. 

Let $G= \SL(2,3) \ltimes E$. Observe that $E/Z$ is a chief factor
of $G$ since $\SL(2,3)$ acts transitively on $(E/Z)^{\#}$.
Let
$\Lambda \in \Irr(E)$ be a character of degree 3. Observe that
$\Lambda$ is $G$-invariant. Since $(|{\bf Q}_8 | , | E |)=1$ 
and $\Lambda$ is $G$-invariant, we have that 
$\Lambda$ extends
to ${\bf Q}_8 \ltimes E$. Since ${\bf Z}_3$ is cyclic and $\Lambda$
is $G$-invariant, 
$\Lambda$ extends to ${\bf Z}_3 \ltimes E$. 
Thus $\Lambda$ extends to $G$
(see Corollary 11.31 of \cite{isaacs}).
Let $\chi\in \Irr(G)$ be an extension of $\Lambda$. We can check 
that $\chi$ is a faithful character. 
By Lemma \ref{itworksifchiextension}
it follows that $\eta(\chi)=1$ and we 
are done with this case.

{\bf (iii) and (iv)} This is Lemma \ref{caseeasier} (ii) and (iii).

{\bf (v)} Using the Atlas, 
we can check that $\Spl(4,3)$ has a subgroup $H$ such that
${\bf F}(H)$ is an
extra-special group of order $2^5$, namely the central product
of  the quaternion group
 ${\bf Q}_8$ and the dihedral group of order 8,
 $|{\bf F}_2(H)/{\bf F}(H)|=5$
and $H/{\bf F}_2(H)$ is a 
subgroup of the cyclic group of order 2. 

Let $E$ be an extraspecial group of order $3^5$ and exponent 3,
and let $\Lambda \in \Irr(E)$ be a character with degree 9.
By Lemma \ref{splitting}, we may assume that $H$ is a subgroup 
of $\Aut(E \modu {\bf Z}(E))$. Set 
$G = H \ltimes E$. Observe that $\Lambda $ is $G$-invariant.
Since $(|H|, |E|)=1$ and $\Lambda$ is $G$-invariant, $\Lambda$
extends to $G$. Let $\chi \in \Irr(G)$ be an extension of
$\Lambda$. Lemma \ref{itworksifchiextension} implies that
 \eqref{basicone} holds for $\chi$ and so we are done.
\end{proof}

\end{subsection}

\end{section}

\begin{section}{Theorem B}
In this section we study the solvable groups $G$ that have an irreducible character
$\chi$ such that $\chi\overline{\chi}$ has exactly 2 distinct non-principal irreducible
constituents. 
Theorem B will be proved in Subsection 4.2. 
Examples of groups satisfying the hypotheses of Theorem B will be provided in 
Subsection 4.2. If, in addition, the group $G$ is
also nilpotent, the structure of the group is very constrained. We will classify
those groups in Subsection 4.3
\begin{subsection}{Hypothesis \ref{hypot}}

\begin{hypo}\label{hypot}

Let $G$ be a finite solvable group.
Let $\chi \in \Irr(G)$ be a faithful character.
   Assume
\begin{equation}\label{basic}
\chi\overline{\chi}= 1_G + m_1\alpha_1 + m_2\alpha_2
\end{equation}

\noindent where $\alpha_1, \alpha_2 \in \Irr(G)$ are non-principal
characters, $m_1$ and $m_2$ are strictly positive integers.
 Set $Z= {\bf Z}(G)$. Let $\lambda \in \Irr(Z)$
be the character such that $[\chi_Z, \lambda]\neq 0$.
\end{hypo}
 \end{subsection}

\begin{subsection}{Derived length}

\begin{lemma}\label{lemma1}
Assume Hypothesis \ref{hypot}. Let $N$ be a normal subgroup of $G$
 Then  $\chi_N \in \Irr(N)$ if and only if 
 $N\not\leq \Ker(\alpha_1)$
and $N \not\leq\Ker(\alpha_2)$.
\end{lemma}

\begin{proof}
 Observe that 
$$[(\chi\overline{\chi})_N, 1_N] = [1_N +m_1(\alpha_1)_N +m_2(\alpha_2)N, 1_N]=
1+m_1[(\alpha_1)_N, 1_N]+ m_2[(\alpha_2)_N, 1_N].$$
Thus $[(\chi\overline{\chi})_N, 1_N]= [\chi_N, \chi_N]=1$ if and only if
$[(\alpha_1)_N, 1_N]=0$ and $[(\alpha_2)_N, 1_N]=0$. Therefore
$[\chi_N, \chi_N]=1$ if and only if  $N\not\leq \Ker(\alpha_1)$
and $N \not\leq\Ker(\alpha_2)$. Since $[\chi_N, \chi_N]=1$ if and only
if
$\chi_N \in \Irr(N)$, the result follows.
\end{proof}

\begin{lemma}\label{thetaregular} 
Let $T$ be a finite  group and $\Theta$ 
a character of $T$. Suppose
that $T=R\times S$. Assume that $\Theta (g) = 0$ if
 $g \in T \setminus R$ or if  $g \in T \setminus S$. Then

\begin{equation}\label{itsregular}
\Theta(g ) = \left\{
\begin{array}{ll} 
                  0 & \mbox{if $g\neq 1$,} \\
               \Theta(1) & \mbox{otherwise.} 
\end{array}\right.
\end{equation}
\noindent Thus $\Theta$ is a multiple of the regular character.
 \end{lemma}
\begin{proof}
Since $\Theta(g)=0$ for $g\in T\setminus R$ and 
$R\cap S =1$,      $\Theta(s)=0$ for
all $s \in S\setminus \{1\}$. Also $\Theta(t)=0$ for all $t \in T\setminus S$.
Thus \eqref{itsregular} holds.
\end{proof}
 
\begin{lemma}\label{kernelcenter}
Assume Hypothesis \ref{hypot}.  
Then either $\Ker(\alpha_1)= Z$ or $\Ker(\alpha_2)= Z$.
\end{lemma}

\begin{proof} 
Assume that $\Ker(\alpha_1) \neq Z$ and 
$\Ker(\alpha_2)\neq Z$. 

	Let $R/ Z $ be a chief factor of $G$, 
where $R \leq\Ker(\alpha_1)$.
Note that such an 
$R$ exists since $Z = \Ker (\alpha_1)\cap \Ker(\alpha_2)$ 
by Lemma \ref{center}, and  $\Ker(\alpha_1) \neq Z$ by assumption. 
Observe also that $ R \not\leq \Ker(\alpha_2)$.
Similarly, let $S/ Z$
 be a chief factor
of $G$ such that $S \leq\Ker(\alpha_2)$. 
Set $T=RS$. Since  
 $\Ker(\alpha_1)\cap \Ker(\alpha_2)= Z$, we have that 
$R\cap S= Z$. Hence 
$T\not\leq \Ker(\alpha_1)$ and $T\not\leq \Ker(\alpha_2)$.
By Lemma \ref{lemma1}, we have that $\chi_T\in \Irr(T)$.

 Since $S/Z$ is a 
chief factor and  $T=RS$ with $ R\cap S=Z$, we have that 
$T/R$ is a chief factor of $G$. 
Let $\psi \in \Irr(R)$ be such that 
$[\chi_R , \psi ]\neq 0$. By Theorem 6.18 of \cite{isaacs}
one of the following holds
	
	a) $\chi_R =\psi$.

	b) $\chi_R = e_1 \psi$ for some $\psi \in \Irr(R)$ and 
some integer $e_1>0$ such that ${e_1}^2=|T:R|$.

	c) $\psi^T=\chi_T$.

Observe that $\chi_R =\psi$ can not hold since $R\leq\Ker(\alpha_1)$
and Lemma \ref{lemma1}.
Observe also that b) and c)  imply that $\chi(g)=0$ if
$g\in T\setminus R$. Similarly, we can check that
$\chi(g)=0$ if $g\in T\setminus S$. Since $R/Z$, $S/Z$ are chief factors of
$G$, they are elementary abelian. Thus we have 
\begin{equation*}
\chi\overline{\chi}(e) = \left\{
\begin{array}{ll} 
                  0 & \mbox{if $e\in T \setminus  R$} \\ 
                 0 & \mbox{if $e\in T \setminus   S$} \\
               \chi(1)^2 & \mbox{ if } e \in R\cap S=Z.
\end{array}\right.
\end{equation*}

            Recall that 
 $R$ and $S$ are normal subgroups of $G$ with 
 $R\cap S= Z$ and $T=RS$. Then we have that
$T/Z= R/Z \times  S/Z$. 
  By  Lemma
\ref{thetaregular} we have that   $\chi\overline{\chi}$ is a multiple of
the 
regular character of $RS \modu Z$, i.e. the regular character 
of $RS/Z$ inflated to $RS$. Thus $\chi(1)^2\geq |T:Z|$. Since
 $\chi_T\in \Irr(T)$
and $Z={\bf Z}(G)$, we have that $\chi(1)^2\leq |T:Z|$. 
Therefore $(\chi\overline{\chi})_T=1_Z^T$ is the regular 
character of $T\modu  Z$.
Because all the characters
of $\Irr(T \modu Z)$ are linear, and thus appear with multiplicity 
1 in $1_Z^T$, and $(\chi\overline{\chi})_T=1_Z^T$, 
we have that $m_1=m_2=1$.
	
Observe that 
\begin{equation*}
\begin{split}
[(\chi\overline{\chi})_R, 1_R]& =
 [((\chi\overline{\chi})_T)_R, 1_R]\\
	&=[(1_Z^T)_R, 1_R]= |T:R|.
\end{split}
\end{equation*}
 Recall that $R\leq\Ker(\alpha_1)$ and $R \not\leq\Ker(\alpha_2)$.
 Applying \eqref{basic}
we get that  $$1+\alpha_1(1) =[(\chi\overline{\chi})_R, 1_R]= |T:R|.$$
 Similarly
we can check that $$1+\alpha_2(1)= |T:S|.$$ Thus	
$$\chi(1)^2= |T:Z|= 1+\alpha_1(1) + \alpha_2(1)=|T:R|+|T:S|-1.$$
Because  $T=RS$ and $R\cap S=Z$, we have that $|T:S|= |R:Z|$. 
Observe that $$|T: Z|=|T:R||R: Z|.$$ 
So we get 
$$|T:R||R: Z|=|T:R|+  |R:Z| -1.$$
\noindent Therefore
$$(|T:R|-1) (|R: Z|-1)= 0.$$
This is impossible since neither $|T:R|$ nor $|R:Z|$ is 1. 
 Thus either 
$\Ker(\alpha_1)=Z$ or $\Ker(\alpha_2)=Z$.
\end{proof}

We are going to describe how $G$ looks, depending
on $\Ker(\alpha_1)$ and $\Ker(\alpha_2)$.

\begin{proposition}\label{kernelsequal}
 Assume Hypothesis \ref{hypot}.
Assume also that $$\Ker(\alpha_1)=\Ker(\alpha_2)=Z.$$
Let $E/Z$ be a chief factor of $G$. 
Then $E/Z$ is a fully ramified section with respect
to $\chi_E$ and $\lambda$. Furthermore
$G/E$ acts faithfully and  symplectically on $E/Z$,
with 2 orbits
on $(E/Z)^{\#}$. Thus $\dl(G)\leq 18$. Also $m_1=m_2=1$,
$\chi(1)^2=|E:Z|$ and $\chi(1)$ is a power of
a prime.
\end{proposition}

\begin{proof} By Lemma \ref{lemma1} we have that
$\chi_E \in \Irr(E)$. Since $\chi_E(1) \neq 1$,  
$E/Z$ is a fully ramified section with respect to 
$\chi_E$ and $\lambda$.  
By Lemma \ref{centralizere} 
we have that ${\bf C}_G(E/{\bf Z}(G))=E{\bf C}_G(E)$ and $E/Z$
has the structure of a symplectic vector space.
Since $\chi$ is a faithful character, we have that 
${\bf Z}(\chi)= {\bf Z}(G)=Z$.
By Lemma \ref{restrictionirreducible} we have that
${\bf C}_G(E) \leq{\bf Z}(\chi)=Z$. Thus ${\bf C}_G(E/Z)=E$.
We conclude
that $G/E$ acts faithfully and symplectically
on $E/Z$. 

	Since $E/Z$ is a fully ramified section with respect to 
$\chi_E$ and $\lambda$, we have that 
\begin{equation}\label{equation2}
(\chi\overline{\chi})_E = 1_E^Z.
\end{equation}
In particular, $\chi(1)^2= |E:Z|$. Since $E/Z$ is a chief factor of
a solvable group, $|E:Z|$ is a power of a prime number.
Therefore $\chi(1)$ is a power of a 
prime number.

	By \eqref{basic} and \eqref{equation2} we have that
\begin{equation}\label{equation1}
1_E^Z= (\chi\overline{\chi})_E = (1_G + m_1\alpha_1 + m_2\alpha_2)_E
=1_E +  m_1(\alpha_1)_E + m_2(\alpha_2)_E.
\end{equation}
By Clifford Theory we have that the irreducible constituents
of $(\alpha_1)_E$ form a $G$-orbit. Similarly, the 
irreducible constituents of $(\alpha_2)_E$ form a
$G$-orbit. If the irreducible constituents of  $(\alpha_1)_E$ are
the irreducible constituents of $(\alpha_2)_E$, then 
by \eqref{equation1}
$G/E$ acts transitively on $1_E^Z \setminus  1_E$. 
Since $\chi_E$ is irreducible, 
by Lemma \ref{itworksifchiextension} and Theorem A 
it follows that 
$\chi\overline{\chi}= 1_G + \alpha$, where $\alpha \in \Irr(G)$.
Thus the irreducible constituents of $(\alpha_1)_E$
and $(\alpha_2)_E$ form two distinct $G$-orbits. 
By \eqref{equation1}, we have that
 the number of orbits of $G$ on $\Irr(E\modu Z)^{\#}$ is 2.
Therefore the number of orbits of $G$ on $(E/Z)^{\#}$ is 2.

Since $G/E$ acts faithfully and irreducibly on $E/Z$,
with 2 orbits on $(E/Z)^{\#}$, by Theorem 
5.4 of \cite{keller} we have that $\dl(G/E)\leq 16$.
 Therefore $\dl(G)\leq 16 +2=18 $.
Observe that $m_1=m_2=1$ by Lemma \ref{a1equal1}.
\end{proof}
	
{\bf Remark.} We will give 
an example of a group satisfying 
the hypotheses of Proposition \ref{kernelsequal}
in Proposition  \ref{examplekernelsequal}.

\begin{proposition}\label{even}
 Assume Hypothesis \ref{hypot}.
Then the order of $G$ is even.
\end{proposition}
\begin{proof} 
Let $\lambda \in \Irr({Z})$ be  such that
$[\chi_{Z}, \lambda]\neq 0$. Let $E/{Z}$ be a chief factor of
$G$.

 If $2 \, \bigm| \, \chi(1)$, then $|G|$ has to be even since 
$\chi(1)\, \big| \,  |G|$. So we may assume that $\chi(1)$ is an odd 
number.
	
	If $\alpha_i$ is a real character, for some $i$,
 then $G$ contains a 
conjugacy class with elements of order 2. 
Thus $G$ has to be of
even order. We may assume that $\alpha_1$ and $\alpha_2$ are not
real characters.

	Since $\chi\overline{\chi}$ is a real character and 
$\alpha_1$ and $\alpha_2$ are not
real characters, by \eqref{basic} we conclude that 
$\overline{\alpha_1}=\alpha_2$. Thus $\Ker(\alpha_1)=\Ker(\alpha_2)$
and $\alpha_1(1)=\alpha_2(1)$.
Since 
${Z}= \Ker(\alpha_1) \cap \Ker(\alpha_2)$ by Lemma \ref{center}, 
we conclude that 
   $\Ker(\alpha_1)=\Ker(\alpha_2)= Z$.
 Thus by the previous  Lemma, $\chi(1)$ is 
a power of an odd prime and $m_1=m_2=1$. 
In particular $4$ divides $\chi(1)^2-1$.

	By \eqref{basic} we have that
\begin{equation*}
	\chi(1)^2= 1 + 2\alpha_1(1).
\end{equation*}
	Thus $2\alpha_1(1)= \chi(1)^2-1$. Since $4$ divides $\chi(1)^2-1$,
 we must have that $2$ divides $\alpha_1(1)$. By Hypothesis \ref{hypot}
we have that $\alpha_1 \in \Irr(G)$. Thus $2$ divides the order of $G$.
\end{proof}
\begin{proposition}\label{abelian}
Assume Hypothesis \ref{hypot}. Assume also that 
$K=\Ker(\alpha_1)$ is an abelian group properly containing 
$Z$. Let $\theta \in \Irr(K)$ be such
 that $[\chi_K, \theta]\neq 0$.
Denote by $ G_{\theta} = \{g \in G | \theta^g=\theta\}$ 
the stabilizer in $G$ of 
$\theta$. Let $L/K$ be a chief factor of $G$.
Then 

(i) $\Ker(\alpha_2)= Z$.

(ii) $G= G_{\theta}L$ and $G_{\theta}\cap L= K$.

(iii) $G_{\theta}/K$ acts faithfully and irreducibly on
$L/K$, and transitively
on $(L/K)^{\#}$.

  In particular,  $\dl(G)\leq 6$ and $\chi(1)$ is a power
of a prime. Also $m_1=1$.

\end{proposition}

\begin{proof}
Since $\Ker(\alpha_1) \neq Z$, by  Lemma \ref{kernelcenter}
it  follows that
$\Ker(\alpha_2)=Z$. 

By Lemma \ref{lemma1} we have that 
$\chi_L \in \Irr(L)$. Since
 $\chi_K$ is reducible,
$K$ is abelian and ${\bf Z}(\chi)=Z\neq K$, by Theorem 6.18 of 
\cite{isaacs} we have that $\chi_L= \theta^L$ and that 
$G_{\theta}< G$.
By Clifford theory it follows that $\chi$ is 
induced from an irreducible character of $G_{\theta}$.
By Exercise 5.7 of \cite{isaacs}, we have that 
$G= G_{\theta}L$. 
Since $L/K$ is an abelian group, $G_{\theta} \geq K$
and $G=G_{\theta}L$, we have that 
 $G_{\theta}\cap L$ 
is a normal subgroup of $G$ containing $K$.
 Observe that 
$\chi_{G_{\theta}}$ is reducible since 
$\chi$ is induced from some character of $G_{\theta}$ and $G_{\theta}<G$.
Thus 
$\chi_{ G_{\theta}\cap L}$ is a reducible 
character, since $ G_{\theta}\cap L\leq G_{\theta}$.
Because $K \leq G_{\theta}\cap L \leq L$,  $G_{\theta}\cap L$ 
is a normal subgroup of $G$ and  $L/K$
is a chief factor of $G$, we have that either 
$G_{\theta}\cap L= L$ or $G_{\theta}\cap L=K$.
If $G_{\theta}\cap L= L$ then $G=G_{\theta}L= G_{\theta}<G$. 
Thus   $G_{\theta}\cap L= K$. 

	Observe that ${\bf C}_G(L/K)$ is a normal subgroup of
$G$ that contains $L$. Set $C= {\bf C}_G(L/K)$.
Since $G = G_{\theta} L$, we have that $C= (C\cap G_{\theta}) L$.
Observe that $C \cap G_{\theta} \trianglelefteq G_{\theta}$ 
and $[C\cap G_{\theta}, L] \leq [C,L] \leq K$. Thus 
$C\cap G_{\theta}\leq G_{\theta}$ 
is a normal subgroup of $G$. 
If $C\neq L$ then 
$C\cap G_{\theta}\leq G_{\theta} $ 
is a normal subgroup of $G$
that properly contains $K$. Thus $\chi_{C\cap G_{\theta}}\in
\Irr(C\cap G_{\theta})$ by Lemma \ref{lemma1}. 
This can not be because $\chi$ is 
induced from some character of $G_{\theta}$ and 
$C\cap G_{\theta} \leq G_{\theta}< G $. Thus 
$C\cap G_{\theta}=K$ and $C=L$. Therefore $G_{\theta}/K$ acts
faithfully on $L/K$. Also $G_{\theta}/K$
 acts irreducibly on $L/K$ since 
$G = G_{\theta} L$ acts irreducibly on the chief factor $L/K$, 
$L=C= {\bf C}_G(L/K)$ and $C\cap G_{\theta}=K$.
 
	By Lemma \ref{basico1} we have that
\begin{equation}\label{structurechi}
(\chi\overline{\chi})_L=1_K^L + \Phi,
\end{equation}
\noindent where $\Phi$ is either the zero function or a character of $L$,
 and $[\Phi_K, 1_K]=0$.
Also, by \eqref{basic} we have that
\begin{equation}\label{structurechi2}
(\chi\overline{\chi})_L=1_L + (m_1\alpha_1)_L + (m_2\alpha_2)_L.
\end{equation}
\noindent Since $\Ker(\alpha_2)=Z < \Ker(\alpha_1)$, we have
that $[(\alpha_2)_K, 1_K]=0$. Thus by \eqref{structurechi} and 
\eqref{structurechi2} it follows that the irreducible 
constituents of $(\alpha_1)_L$ form the set
$\Irr(L \modu K)^{\#}$. It follows that 
 $G$ acts transitively on 
$\Irr(L \modu K)^{\#}$. 
Therefore $G$ acts transitively on $(L/K)^{\#}$. 

Since ${\bf C}_G(L/K)=L$, we have that 
that $G/L$ acts faithfully on $L/K$ and transitively on $(L/K)^{\#}$.
Using  Lemma \ref{huppert}, we can check that $\dl(G/L) \leq 4$.  
Since $L/K$ is a chief factor of $G$ and $K$ is abelian, it follows
that $\dl(L)\leq 2$. We conclude that 
$\dl(G) \leq 6$.
	
Since $\chi_L= \theta^L$ and $\theta(1)=1$, we have that
 $\chi(1)=|L:K|$. Observe that $|L:K|$ is a power of a prime 
since 
$L/K$ is a chief factor 
of a solvable group. 
Thus $\chi(1)$ is a power of a prime number.
By Lemma \ref{a1equal1} we have that $m_1=1$. 
\end{proof}

{\bf Remark.} In Proposition \ref{examplesabeliancase},
we will prove that for any prime $p$
and any nonzero positive integer $m$,  there exist
a group $G$ and a character $\chi \in \Irr(G)$ with $\chi(1)=p^m$
  satisfying the hypotheses of Proposition \ref{abelian}.

\begin{lemma}\label{case9}
Assume that $G$ is   a group such that ${\bf F}(G)$ is an
extra-special group of order $2^5$, $|{\bf F}_2(G)/{\bf F}(G)|=5$
and $G/{\bf F}_2(G)$ is a subgroup of the cyclic group of order 2.

Then $\cd(G) \subseteq  \{1, 2, 4, 5, 8, 10 \}$.
Also, if $\alpha \in \Irr(G)$ and $\alpha (1)=8$, then 
$\Ker(\alpha)=1$.

\end{lemma}
\begin{proof}
Observe that ${\bf Z}(G)= {\bf Z}({\bf F}(G))$. 
Since  ${\bf F}(G)$ is an
extra-special group of order $2^5$, there exists a unique
$\Lambda \in \Irr({\bf F}(G))$ such that $\Lambda(1)=4$. This 
$\Lambda$ is a $G$-invariant character. Since the 
$p$-Sylow subgroups of $G/ {\bf F}(G)$ are cyclic 
for all primes $p$, and 
$\Lambda$ is a $G$-invariant character, we have that $\Lambda$
extends to $G$. Observe that $\cd(G/ {\bf F}(G)) \subseteq \{1,2\}$.
Thus the degrees of the characters in $\Irr( G | \Lambda)$
lie in the set  $\{4,8\}$, while those of $\Irr( G \modu {\bf F}(G))$ 
lie in $\{1,2\}$.

Let $\lambda \in \Lin ({\bf F}(G)\modu {\bf Z}(G))^{\#}$ be a linear 
character of ${\bf F}(G)$.  Since ${\bf F}_2(G)/{\bf F}(G)$ 
acts without fixed points on $({\bf F}(G)/{\bf Z}({\bf F}(G)))^{\#}$,
i.e on $\Lin ({\bf F}(G))^{\#}$, we have that $\lambda$ induces 
irreducibly to ${\bf F}_2(G)$. 
If $\chi \in \Irr(G)$ lies above $\lambda$, then $\chi_ {{\bf F}(G)}$
is the sum of $\chi(1)$ distinct linear characters of
 ${\bf F}(G)\modu {\bf Z}(G)$. Because of $|G:{\bf F}(G)| \in \{ 5, 10\}$,
 we have that 
$\chi(1)\in \{ 5, 10 \}$.
\end{proof}
\begin{proposition}\label{knotabelian}
Assume Hypothesis \ref{hypot}. Assume also that $\Ker(\alpha_1)$ is not 
abelian.
Let $E/Z$ be a chief factor of $G$ such that $E\leq \Ker(\alpha_1)$. 
Let 
$F_2/E$ be  the Fitting subgroup of  $G/E$.
Then $m_1=1$, $m_2 =4$, $\alpha_1(1)=3$ and $F_2/E$
is isomorphic to
$Q_8$. Also, one of the following holds
 
(i) $\chi(1)=10$, $G/E \cong \SL(2,3)$.

(ii) $\chi(1)=14$, $G/E \cong \widetilde{\GL}(2,3)$.
 
\end{proposition}
\begin{proof}
Set $K=\Ker(\alpha_1)$.
Let $\theta \in \Irr(K)$ be a character such that 
$[\chi_{K}, \theta] \neq 0$.
Since $K$ is nonabelian and $\chi$ is faithful, $\chi_K$ is a sum of 
$G$-conjugate nonlinear 
characters, one of which is $\theta$. 
Therefore $\theta(1)>1$.

Let $\lambda \in \Irr(Z)$
be such that $[\lambda, \theta_Z]\neq 0$.
By Lemma \ref{kernelcenter} we have that $\Ker(\alpha_2)= Z$.
 
\begin{claim}\label{maximal}
For any normal subgroup $N$ of $G$ such that $Z< N \leq K$, we have 
$\theta_N \in \Irr(E)$. In particular
\begin{equation}\label{thetarestrictsirreducibly}
 \Lambda =\theta_E \in \Irr(E).
\end{equation}
\end{claim}
\begin{proof}
Suppose that  $\theta_N$ is a reducible character. Then 
$[\theta_N, \theta_N] >1$.
Since $[\chi_K, \theta]\neq 0$, we have
that  
$$[(\chi\overline{\chi})_N, 1_N]=[\chi_N, \chi_N] > [\chi_K, \chi_K]
=[(\chi\overline{\chi})_K, 1_K].$$ 
Thus there exists an irreducible constituent
$\alpha$  of $\chi\overline{\chi}$ 
such that $N\leq \Ker(\alpha)$ and $K \not\leq \Ker(\alpha)$.
But the irreducible constituents of $\chi\overline{\chi}$ are
$1_G$, $\alpha_1$ and $\alpha_2$ and their kernels are
$G$, $K$ and $Z$ respectively. We conclude that $\theta_N$ is an irreducible 
character.
\end{proof}

 Since 
$\theta_E$ is not a linear character, $E/Z$ is a chief factor 
of $G$ and $Z={\bf Z}(G)$, the section $E/Z$ is  
fully ramified with respect
to $\Lambda$ and $\lambda$. Thus 
\begin{equation}\label{lambdaisginvariant}
\mbox{$\Lambda$ is $G$-invariant and $\Lambda(1)^2=|E:Z|$}.
\end{equation}
\begin{claim}\label{faithful}
 ${\bf C}_G(E/Z)=E$. Thus $G/E$ acts faithfully 
 on the symplectic vector space $E/Z$. Also $G/E$ acts 
 transitively on $(\Lin(E/Z))^{\#}$, hence on
$(E/Z)^{\#}$.
\end{claim}
\begin{proof}
We will first prove that ${\bf C}_G(E)=Z$. Observe that $C_G(E)$ is 
a normal subgroup of $G$.
Suppose that ${\bf C}_G(E)\neq Z$.
If ${\bf C}_G(E)\leq K$, by Claim  \ref{maximal}  we have that 
$\theta_{{\bf C}_G(E)}\in \Irr(C_G(E))$.
By Lemma \ref{restrictionirreducible} we have that
$$E\leq {\bf C}_K ({\bf C}_G(E))\leq {\bf Z}(\theta).$$
Since  $\Lambda =\theta_E \in \Irr(E)$ and $E\leq{\bf Z}(\theta)$, it
follows that $\Lambda$ is a linear character. Therefore $\theta(1)=1$.
But $\theta(1)>1$. 
So we may assume that ${\bf C}_G(E)\not\leq K$.
 By Lemma \ref{lemma1} we have that 
$\chi_{{\bf C}_G(E)}\in \Irr({\bf C}_G(E))$.
As before, by Lemma \ref{restrictionirreducible} we have that 
$$E\leq {\bf C}_G ({\bf C}_G(E))\leq {\bf Z}(\chi)=Z<E,$$
\noindent a contradiction. Thus ${\bf C}_G(E)=Z$.
By Lemma \ref{centralizere} we have that ${\bf C}_G(E/Z)={\bf C}_G(E) E$.
Since ${\bf C}_G(E)=Z$, it follows that  ${\bf C}_G(E/Z)=E$.

By \eqref{thetarestrictsirreducibly} and \eqref{lambdaisginvariant}
we have that $E/Z$ is a fully ramified section with respect to 
$\Lambda$ and $\lambda$. Since $Z= {\bf Z}(G)$, we have 
\begin{equation}\label{lambdarestrictsproduct}
 (\theta \overline{\theta})_E= \Lambda \overline{\Lambda}=1_Z^E.
\end{equation}
Since $[\Lambda, (\chi)_E]\neq 0$, we have 
$(\chi\overline{\chi})_E= 1_Z^E + \Phi,$ 
where $\Phi$ is character of $E$. Since
$K=\Ker(\alpha_1)\geq E$ and \eqref{basic} holds we have
\begin{equation*}
\begin{split}
(\chi\overline{\chi})_E &= (1_G+ m_1\alpha_1 +m_2\alpha_2)_E\\
			&= (1+m_1 \alpha_1(1))1_E + m_2 (\alpha_2)_E\\	
&= 1_Z^E + \Phi.
\end{split}
\end{equation*}
By Clifford Theory, the irreducible constituents of $(\alpha_2)_E$
are $G$-conjugates. Thus $(\Irr(E \modu Z))^{\#}$ is  a
$G$-orbit.
\end{proof}

\begin{claim}\label{symplecticoptions2}
Set $e^2=|E:Z|$ and $\bar{G}= G/E$. 
Then one of the following holds:
 
(i) $e=2$, and either $\bar{G}\cong S_3$ or $\bar{G}\cong {\bf Z}_3$.

	(ii) $e=3$, and 
either $\bar{G}={\bf Q}_8$ or $\bar{G}\cong \SL(2,3)$.

         (iii) $e=5$ and
$\bar{G}\cong \SL(2,3)$.

	(iv) $e=7$ and 
$\bar{G}\cong \widetilde{\GL}(2,3)$.

	(v) $e=9$, ${\bf F}(\bar{G})$ is an
 extra-special group of order $2^5$, 
 $|{\bf F}_2(\bar{G})/{\bf F}(\bar{G})|=5$
 and $\bar{G}/{\bf F}_2(\bar{G})$ is a 
 subgroup of the cyclic group of order 2.

\end{claim}
\begin{proof} This follows from Claim \ref{faithful} and 
Theorem  \ref{symplectictrans}.
\end{proof}

\begin{claim}\label{lambdaextends}
Let $e \in \{2,3,5,7\}$.
 Set $\Lambda =(\theta)_E$. Then $\Lambda\in \Irr(E)$, $\Lambda(1)=e$ 
and $\Lambda$ extends to $G$, i.e. there exists $\Lambda^e \in \Irr(G)$
such that $(\Lambda^e)_E=\Lambda$.
\end{claim}

\begin{proof}
By \eqref{thetarestrictsirreducibly} we have that 
$\Lambda= \theta_E \in \Irr(E)$. By \eqref{lambdarestrictsproduct}
we have that $\Lambda(1)=e$.

Assume $e=2$. Since the $p$-Sylow subgroups of $G/E$ are cyclic
for all  
primes $p$, by Corollary 11.22 and Corollary 11.31 of 
\cite{isaacs} we have that  $\Lambda$ extends to $G$. 

Assume now that $e \in \{3,5,7\}$. Observe that for primes $p\neq 2$,
the $p$-Sylow subgroups of $G/E$ are cyclic.
If we can prove that $\Lambda$ extends to $ET$, where
$T$ is a 2-Sylow subgroup of $G$, then  $\Lambda$
would extend to $G$ (see Corollary 11.31 of \cite{isaacs}). 
In our situation 
we have that either
$TE/E\cong{\bf Q}_8$ or 
$TE/E\cong{\bf Q}_{16}$.
It follows from the example on page 301 of
\cite{reiner} 
that $H^2({\bf Q}_8, {\bf C})=1$ and  $H^2({\bf Q}_{16}, {\bf C})=1$.
By Proposition 11.45 of \cite{reiner} this implies 
that $\Lambda$ extends
to $ET$. Thus $\Lambda$ extends to $G$.
\end{proof}

\begin{claim}\label{psiisthere}
Let $e \in \{2,3,5,7\}$.
Let  $\Lambda^e \in \Irr(G)$ be an extension
 of $\Lambda\in \Irr(E)$. Then
\begin{equation}\label{chiproduct}
\chi= \Lambda^e \psi
\end{equation}
\noindent for some $\psi \in \Irr(G \modu E)$. Also

	(i) $\Lambda^e \overline{\Lambda^e}= 1_G + \alpha_2$.

	(ii)  $\psi \overline{\psi}= 1_G + \alpha_1$.

        (iii) $\alpha_1 \alpha_2 = \alpha_1(1) \alpha_2$ and 
$m_2=1+\alpha_1(1)$.
\end{claim}
\begin{proof}
Since $[\Lambda, \chi_E]\neq 0$ and $\Lambda$ extends to 
$\Lambda^e \in \Irr(G)$, by Corollary
6.17 of \cite{isaacs} there exists $\psi \in \Irr(G\modu E)$
such that \eqref{chiproduct} holds. 

Observe that $\Ker(\Lambda^e)\leq E$ since $\Lambda^e$ is an extension
of $\Lambda$ and $G/E$ acts faithfully on $E/Z$ (see Claim \ref{faithful}). 
Since $\chi$ is faithful and $\Ker(\Lambda^e)\leq E$, by
\eqref{chiproduct} we have that $\Ker(\Lambda^e)=\Ker(\chi) \cap E =1$.
Thus $\Lambda^e$ is a faithful character of $G$, $E/Z$ is a chief factor of 
$G$, where $Z={\bf Z}(G)$,  and $G/E$ acts on 
$(E/Z)^{\#}$ transitively (see Claim \ref{faithful}).
By Lemma \ref{itworksifchiextension}, we conclude that
	$$\Lambda^e \overline{\Lambda^e}= 1_G+ c_1 \gamma$$
\noindent for some  $\gamma \in \Irr(G)$ and some integer $c_1$. 
By Theorem A we have that $c_1=1$.
 By Lemma \ref{center}, we observe
that 
\begin{equation}\label{gammakernelz}
\Ker(\gamma)=Z. 
\end{equation}

	Let 
$\psi \overline{\psi} = 1_G + \sum_{i=1}^{\eta(\psi)} d_i\delta_i$,
for some integers $d_i>0$  and some 
distinct characters $\delta_i \in \Irr(G)$ for $i=1, \ldots , \eta(\psi)$. 
Since $\psi \in \Irr(G \modu E)$, we have that 
$\delta_i \in \Irr(G \modu E)$ for all $i$.
Observe 
\begin{equation*}
 \begin{split}
\chi \overline{\chi} &= \Lambda^e \psi  \overline{\Lambda^e \psi} \\
  & = (\Lambda^e \overline{\Lambda^e}) (\psi \overline{\psi})\\
 & = (1_G + \gamma)(  1_G + \sum_{i=1}^{\eta({\psi})} d_i\delta_i)\\
 & = 1_G +\gamma +  \sum_{i=1}^{\eta({\psi})} d_i\delta_i +  
\sum_{i=1}^{\eta({\psi})} d_i \gamma \delta_i\\
&  = 1_G+ m_1 \alpha_1 + m_2 \alpha_2
\end{split}
\end{equation*}
\noindent where the last equality is \eqref{basic}. Since 
$\Ker(\delta_1)\geq E$
and $\Ker(\gamma) =Z \not\leq E$ by \eqref{gammakernelz}, we have that 
$\gamma \neq \delta_1$. It follows that $\eta(\psi)=1$.  

Since $\Ker(\delta_1)\geq E$ and  
$\Ker(\gamma)=Z<E$,
we have that $\Ker(\gamma \delta_1) \not\geq E$.
Since $m_1=1$, $\Ker(\alpha_1)\geq E$ and $\Ker(\alpha_2)=Z$,
it 
follows that $\delta_1=\alpha_1$, $\gamma=\alpha_2$ and $\gamma \delta$ is 
a multiple of $\alpha_2$. Thus $\gamma\delta = \alpha_1 (1) \alpha_2$.
Therefore $m_2= 1 +\alpha_1(1)$. 
\end{proof}

\begin{claim}\label{howg/elooks}
Let 
$e \in \{2,3,5,7\}$.
 Set
 $H=\Ker(\psi)$, where $\psi \in \Irr(G \modu E)$ is as in Claim 
\ref{psiisthere}. Let $Z_H/ H = {\bf Z}(G/H)$.
 Then there exists $L\trianglelefteq G$ such that
 $L/Z_H$ is  a chief factor of $G$.
 Also $E< L$,  $\psi(1)^2 =|L:Z_H|$ and one of the following holds

(i) $\psi(1)=2$ and either ${G/L}\cong {\bf Z}_3$
or $G/L \cong S_3$. 

(ii) $\psi (1)=3 $ and either 
${G/L}\cong  {\bf Q}_8$ or ${G/L} \cong \SL(2,3)$.

(iii) $\psi(1)=5$ and  ${G/L}\cong \SL(2,3)$.
 
(iv) $\psi (1)=7$ and 
${G/L} \cong \widetilde{\GL}(2,3)$.

(v) $\psi(1)= 9$, ${\bf F}(G/L)$ is an
extra-special group of order $2^5$, $|{\bf F}_2(G/L)/{\bf F}(G/L)|=5$
and $(G/E)/{\bf F}_2(G/L)$ is a 
subgroup of the cyclic group of order 2.
\end{claim}
\begin{proof}
By Claim \ref{psiisthere} (ii), we have that $\psi(1)^2 = 1 +\alpha(1) >1$.
Thus $\psi(1)>1$. Since $\psi$ is a faithful irreducible
character of $G/H$ and $\psi(1)>1$, it follows that $G/H$ is not
abelian. Therefore there exists a normal subgroup $L$ of $G$ 
such that $L/Z_H$ is a chief factor of $G$.
 
Observe that $\psi$ is a faithful and irreducible character of 
$G/\Ker(\psi)=G/H$. By Claim \ref{psiisthere} (ii), it follows
that $\eta(\psi)=1$.  Thus Theorem A implies 
Claim  \ref{howg/elooks}.
\end{proof}

\begin{claim} $e\neq 2$.
\end{claim}

\begin{proof}
	Assume that $e=2$. Then $G/E \cong S_3$ or $G/E \cong {\bf Z}_3$ by 
Claim \ref{symplecticoptions2}.
By Claim \ref{howg/elooks} we have that there exists $L\triangleleft G$
such that $E<L$ and $G/L$ is of the form  (i), (ii), (iii), (iv) or (v).
In this case, $G/L\cong {\bf Z}_2$ or $G=L$ since the only non-trivial 
normal subgroup of $G/E\cong S_3$ or $G/E\cong {\bf Z}_3$  is
 ${\bf Z}_3$. And ${\bf Z}_2$ is not 
among the possibilities given by Claim \ref{howg/elooks}. So $e \neq 2$.
\end{proof}

\begin{claim} $e\neq 3$.
\end{claim}
\begin{proof}
Assume that $e=3$. By Claim \ref{symplecticoptions2} we have that 
${G/E}\cong  {\bf Q}_8$ or 
${G/E} \cong \SL(2,3)$.  
By Claim \ref{howg/elooks}, we have that $E<L$, where  $L$ is a
normal subgroup of $G$ such that  ${\bf F}(G/L)$ is isomorphic to 
either ${\bf Z}_3$
or ${\bf Q}_8$.
This is impossible if $G/E \cong {\bf Q}_8$. Thus
  $G/E \cong \SL(2,3)$, $G/L\cong {\bf Z}_3$ and
$\psi(1)=2$. Let $M$ be the  normal subgroup of $G$ such that
$E< M$ and $M/E={\bf F}(G/E) \cong Q_8$.
By Claim \ref{psiisthere} (ii) we have that
$\alpha_1(1)=3$ and $\alpha_1 \in \Irr(G \modu E)$.
 Since  $\alpha_1 \in \Irr(G \modu E)$, $G/E \cong \SL(2,3)$ and
$|G/M| = 3$, we have that   $\alpha_1$ is induced from a linear
character of $M$. In particular 
\begin{equation}\label{alpha1zerooutside}
\alpha_1(g)=0 \mbox{ for any }
g \in G \setminus M.
\end{equation}
	Since $e=3$, by Claim \ref{lambdaextends} we have that 
 $\Lambda^e(1)=3$. By Claim \ref{psiisthere} (i) we have that
$\alpha_2(1)=8$. Since $M$ is a normal 
subgroup of $G$, $(   |G/M|, \alpha_2(1)  ) =1$, and 
$\alpha_2 \in \Irr(G)$, by exercise 6.7 of \cite{isaacs} we 
have  
 that $(\alpha_2)_M  \in \Irr(M)$. 
Thus  $(\alpha_2)_M \in \Irr(M \modu Z)$,
$(\alpha_2)_M(1)=8$,
$E/Z$ is elementary abelian of order $9$, $|M/E|=8$ and $E$ is a normal
subgroup of $G$.  By Clifford Theory we conclude  that $(\alpha_2)_M$
is induced from a linear character of $E$. Therefore
\begin{equation} \label{alpha2zerooutside}
\alpha_2(g) =0 \mbox{  for any } g \in M \setminus E. 
\end{equation}
Since $E$ is normal in $G$, while  $(\alpha_2)_M\in \Irr(M)$ is induced
from a linear character of $E$ and $\alpha_2(1)=8$,
we have that $(\alpha_2)_E$ is the sum
of 8 distinct linear characters of $E$. 
Since $\alpha_2 \in \Irr(G \modu Z)$
and $|E/Z|=9$ and $(\alpha_2)_E$ is the sum
of 8 distinct linear characters of $E$, 
it follows that $(\alpha_2)_E= 1_Z^E -1_E$. Therefore
\begin{equation}\label{valuealpha2one}
\alpha_2 (g) = -1 \mbox{ for any } g \in E \setminus Z.
\end{equation}
By Claim \ref{psiisthere} (iii)
we must have that $\alpha_1 \alpha_2= 3 \alpha_2$.
By \eqref{alpha1zerooutside} and \eqref{alpha2zerooutside}
we have that $\alpha_1\alpha_2 (g) =0$ for any  $g \in G \setminus E$.
Thus 
\begin{equation}\label{alpha2zerog/e}
\alpha_2 (g) =0 \mbox{ for any } g \in G \setminus E.
\end{equation}
Therefore
\begin{equation*}
 \begin{split}
[\alpha_2, \alpha_2] & = \frac{1}{|G|} 
 \sum_{g\in G} \alpha_2(g) \overline{\alpha_2(g)}\\
 & =  \frac{1}{|G|} 
 \sum_{g\in E} \alpha_2(g) \overline{\alpha_2(g)}+  
 \frac{1}{|G|}\sum_{g\in G \setminus E} \alpha_2(g) \overline{\alpha_2(g)}\\
 & =  \frac{1}{|G|} 
 \sum_{g\in E} \alpha_2(g) \overline{\alpha_2(g)} \ 
\text{ by \eqref{alpha2zerog/e}} \\
 & = \frac{1}{|G|} \sum_{g\in Z} \alpha_2(g) \overline{\alpha_2(g)}+
 \frac{1}{|G|} \sum_{g\in E\setminus Z} \alpha_2(g) \overline{\alpha_2(g)}\\
 &= \frac{1}{|G|} |Z| \alpha_2(1)^2 +  \frac{1}{|G|}(|E|-|Z|)  \ 
 \text{\  because $\alpha_2 \in \Irr(G \modu Z)$
 and \eqref{valuealpha2one} }\\
 &= \frac{8^2}{|G/Z|} + \frac{(9-1)}{|G/Z|} \text{ \ since $\alpha_2(1)=8$ and
 $|E/Z|=9$ }\\
 & =\frac{8 \times 9}{24\times 9}=\frac{1}{3}. 
\end{split}
\end{equation*} 

 Since  $\alpha_2 \in \Irr(G)$,
we have that $[\alpha_2, \alpha_2]=1$. But
$[\alpha_2, \alpha_2]= \frac{1}{3} \neq 1$. Thus this
case can not arise and $e \neq 3$.
\end{proof}
\begin{claim} If $e=5$ then $\psi(1)=2$, $\alpha_1(1)=3$, $\alpha_2(1)=24$,
$m_2=4$ and $G/E \cong \SL(2,3)$. This gives case (i) in Proposition
\ref{knotabelian}.
\end{claim}
\begin{proof}
By Claim \ref{symplecticoptions2} we have that $G/E\cong \SL(2,3)$.
Analyzing the possibilities given by Claim 
\ref{howg/elooks}, we have that 
$\psi(1)=2$.
Thus  by Claim \ref{psiisthere} (ii) we have that 
$\alpha_1(1)=3$. Similarly, since $\Lambda^e(1)= 5$, it follows 
by Claim \ref{psiisthere} (i) that
$\alpha_2(1)=24$ and 
by \eqref{chiproduct} that 
$\chi(1)= \psi(1) \Lambda^e(1)=10$.
 By Claim \ref{psiisthere} (iii) we have that
$m_2= 1+ \alpha_1(1)=1+3 =4$ and we are done.
\end{proof}

\begin{claim} If $e=7$ then $\psi(1)=2$, $\alpha_1(1)=3$, $\alpha_2(1)=48$,
$m_2=4$ and $G/E \cong \widetilde{\GL}(2,3)$. This gives case (ii) 
in Proposition \ref{knotabelian}.
\end{claim}

\begin{proof}
As in the proof of the previous claim, we have that $\psi(1)=2$ and
$\alpha_1(1)=3$. Since 
 $\Lambda^e(1)= 7$, by  Claim \ref{psiisthere} (i)  we have that
$\alpha_2(1)=48$.  Since $\Lambda^e(1)= \Lambda(1)=7$, by Claim
\ref{psiisthere} (ii) we have that $\alpha_2(1)= 48$ and 
by \eqref{chiproduct} that 
$\chi(1)= \psi(1) \Lambda^e(1)=14$.
 By Claim \ref{psiisthere} (iii) we have that
$m_2= 1+ \alpha_1(1)=1+3 =4$.
\end{proof}

It remains to analyze the case when $e=9$. We will use another approach 
to that case. The reason is 
that  $H^2( E, {\bf C})\neq 1$ for an extra-special group $E$
of order $2^5$. Thus the character $\Lambda$ may not extend to $G$.
 We start with some statements about  
$\chi$, $\alpha_1$ and $\alpha_2$. Then we use that information 
to prove that $e=9$ can not hold.

 Since $E\leq K$ and $K=\Ker(\alpha_1)$, the restriction $\chi_E$ is
reducible. Since $E/Z$ is a fully ramified section of $G$ with respect
to $\theta_E$, we have that 

	\begin{equation}\label{sizechi}
\chi(1)=em\ \ \mbox{ where}\   m\geq 2.
\end{equation}

Since $E/Z$ is a fully ramified section of $G$ with respect
to $\theta_E$ and $Z={\bf Z}(G)$, $\chi(g)=0$ for all 
$g \in E \setminus Z$.
Thus  $\chi\overline{\chi}(g)=0=1+\alpha_1(1)+m_2\alpha_2(g)$
for all $g \in E \setminus Z$.
Therefore for all $g \in E \setminus Z$ we have
\begin{equation}\label{alpha2termsalpha1}
\alpha_2(g) =- \frac{1+\alpha_1(1)}{m_2}.
\end{equation}

Recall that  
$\alpha_2 \in \Irr(G)$. So $\alpha_2(g)$ is 
an algebraic integer for any $g\in G$.
Since $\frac{1+\alpha_1(1)}{m_2}$ is 
a rational  algebraic integer, it is an integer.
Thus we conclude that
\begin{equation}\label{a2divides}
	m_2 \ \mbox{ divides}\ 1+\alpha_1(1).
\end{equation}

By \eqref{basic}  we have that 
\begin{equation}\label{alpha2intermsalpha1andchi}
\alpha_2(1)= \frac{\chi(1)^2}{m_2} - \frac{1+\alpha_1(1)}{m_2}.
\end{equation}

By \eqref{alpha2termsalpha1} and \eqref{alpha2intermsalpha1andchi} we 
have that 
\begin{equation*}
(\alpha_2)_E= s \, 1_Z^E - \frac{1+\alpha_1(1)}{m_2} 1_E 
\end{equation*}
\noindent for some integer $s \geq \frac{1+\alpha_1(1)}{m_2}$. So 
\begin{equation*}
s e^2 - \frac{1+\alpha_1(1)}{m_2}= 
\alpha_2(1) =\frac{ m^2 e^2}{m_2} -\frac{1+\alpha_1(1)}{m_2}.
\end{equation*}  
\noindent Hence $s = m^2/ m_2$. Therefore 
\begin{equation}\label{a2divideschi} 
m_2 \mbox{ divides } m^2.
\end{equation}
Since $\alpha_2 \in \Irr(G)$ and $\Ker(\alpha_2)=Z$, we have that

\begin{equation}\label{alpha2}
 \alpha_2(1)^2<|G:Z| \ \mbox{ and} \ 
 \alpha_2(1) \ \mbox{ divides }\  |G:Z|.
\end{equation}
   
\begin{claim}\label{EneqK} Assume $e=9$. Then 
$E< K=\Ker(\alpha_1)$. In particular, $\alpha_1$ is not a faithful
character of $G\modu E$.
\end{claim}
\begin{proof}
Suppose $E=K$. 
Observe that ${\bf Z} ({\bf F}(\bar{G}))$ is a non-trivial
 cyclic
group by Claim \ref{symplecticoptions2} (v). 
Therefore there exists a cyclic 
chief factor $H/E$ of $G$ with $H/E \leq {\bf Z} ({\bf F}(\bar{G}))$.
Since $\Lambda \in \Irr(E)$ is $G$-invariant and $H/E$ is 
cyclic, $\Lambda$ extends to $H$. We are assuming that
$E=K=\Ker(\alpha_1)$. Thus $H\not\leq K$,
 and $H\not\leq \Ker(\alpha_2)$. Hence
$\chi_H\in \Irr(H)$ by Lemma \ref{lemma1}. Since $\Lambda$ extends to $H$,
 $H/Z$ is
cyclic and $[\chi_E, \Lambda]\neq 0$, it follows that  
 $\chi_E=\Lambda$. By Claim \ref{lambdaextends} we have that 
 $\Lambda(1)=e$.
Thus  $\chi(1)=e$. But that is  a contradiction with 
\eqref{sizechi}. 
\end{proof}
\begin{claim} $e \neq 9$.
\end{claim}
\begin{proof}
Assume that $e=9$. By Claim \ref{symplecticoptions2} (v) we have that
\begin{equation}\label{indexzg}
|G:Z| \mbox{ divides }  10\times 32 \times 81=25920.
\end{equation}
 By Claim \ref{faithful}, Claim \ref{symplecticoptions2} and 
Lemma 
\ref{case9} we have that $\alpha_1(1) \in \{1,2,4,5,8,10\}$.
By Claim \ref{EneqK} and Lemma \ref{case9}
 we have that $\alpha_1(1) \neq 8$.
Thus
\begin{equation}\label{optionsalpha1}
\alpha_1(1) \in \{1,2,4,5, 10\}.
\end{equation}
This and \eqref{a2divides} imply that 
\begin{equation}\label{a2constraint}
  m_2 \in\{1,2,3,5,6,11 \}.
\end{equation}
Therefore 
\begin{equation}\label{bounding9}
	\frac{1+\alpha_1(1)}{m_2}\leq 11.
\end{equation}
	
	If $m_2 = 1$, then by \eqref{sizechi} and 
\eqref{alpha2intermsalpha1andchi}
 we have
that
	\begin{equation*}
\begin{split}
 \alpha_2(1)& =\frac{\chi(1)^2}{m_2}- \frac{1+\alpha_1(1)}{m_2}\\ 	
	    & \geq    4 \times  81 - 11 =313.
\end{split}
\end{equation*}
But then 
$\alpha_2 (1)^2\geq 313^2= 97969>25920 \geq |G:Z|$.
This is impossible by \eqref{alpha2}. So $m_2 \neq 1$.

	If $m_2 = 2$, then by \eqref{a2divides} we have that 
$\alpha_1(1)=1$ or $\alpha_1(1)=5$.
By \eqref{a2divideschi} we have that  $2 \mid m$
If $m>2$ then $m\geq 4$ and 
 $$\frac{\chi(1)^2}{m_2}= \frac{m^2 9^2}{2}\geq
\frac{ 4^2 \times 9^2}{2}=  8\times 81.$$ 
Thus by \eqref{alpha2intermsalpha1andchi}
we have that 
$\alpha_2(1)\geq 8\times 81 - 3$ and therefore
$\alpha_2(1)^2>|G:Z|$. So we may assume that $m=2$.
Since $\alpha_1(1)\in \{1, 5\}$ and $m_2=2$,
 by \eqref{alpha2intermsalpha1andchi}
  either $\alpha_2(1)= 2\times81-1=161$ or
$\alpha_2(1)= 2\times81-3=159$.
But neither $161$ nor $159$ divides $|G:Z|$.
Hence $m_2 \neq 2$.

	Suppose $m_2 \in \{3, 5, 6, 11\}$. 
By \eqref{a2divideschi} we have that $m_2 \bigm| m^2$. Thus $m\geq 9$.
 Therefore 
 $$\frac{\chi(1)^2}{m_2}\geq  \frac{9^2 9^2}{3}= 2187 .$$
Thus $\alpha_2(1)\geq 2187-11$
by  \eqref{bounding9} and \eqref{alpha2intermsalpha1andchi}.
By \eqref{indexzg} $\alpha_2(1)^2 = 2176^2 >25920\geq |G:Z|$,
a contradiction with \eqref{alpha2}. We conclude that
$m_2\not\in \{3,5,6,11\}$.

We have considered all the possible options given by \eqref{a2constraint} and
we have showed that none of them is possible. We conclude that $e\neq 9$. 
\end{proof}

We have analyzed all the options given
by  Claim \ref{symplecticoptions2}, so the proof is complete.
\end{proof}

We summarize some of the information we obtained in the previous lemmas.

\begin{teoremab}
Assume Hypothesis \ref{hypot}. Then 

(i) The order of $G$ is even.

(ii) $\dl(G) \leq 18$.

(iii) Either $\Ker(\alpha_1)=Z$ or
$\Ker(\alpha_2)=Z$.
 If $\, \Ker(\alpha_1)$ and $\Ker (\alpha_2)$ are both abelian subgroups, then
$\chi(1)$ is a power of a prime number. Otherwise
$\{m_1,m_2\}=\{1,4\}$ and $\chi(1)\in \{10,14\}$. 
\end{teoremab}

\begin{proof} The order of $G$ is even by
Proposition \ref{even}. 
Either $\Ker(\alpha_1)=Z$ or $\Ker(\alpha_2)=Z$ by 
Lemma \ref{kernelcenter}.  Propositions \ref{kernelsequal}, \ref{abelian} and 
 \ref{knotabelian}
imply (ii) and (iii).
\end{proof}

\end{subsection}

\begin{subsection}{Examples}
In this subsection, we study  examples of groups satisfying Hypothesis
\ref{hypot}. 

\begin{lemma}\label{basicequalkernel}
Let $G$ be a finite group and $E$ be a nonabelian 
subgroup. 
Assume   $E/{\bf Z}(G)$ is an abelian 
chief factor of $G$ 
and  $G/E$ acts on $(E/{\bf Z}(G))^{\#}$ 
with 2 orbits. Let $\chi \in \Irr(G)$ be a faithful character.
Assume also that $\chi_E \in \Irr(E)$. Then 
$$\chi\overline{\chi}= 1_G + \alpha_1 + \alpha_2$$
\noindent where $\alpha_1$ and $\alpha_2$ are distinct
irreducible
characters of $G$ and $\Ker(\alpha_1)=\Ker(\alpha_2)={\bf Z}(G)$.
\end{lemma}
\begin{proof}
 Set $Z={\bf Z}(G)$. Let $\lambda \in \Irr(Z)$ be 
such that
$[\lambda, \chi_Z]\neq 0$. Since  
 $\chi_E \in \Irr(E)$ is a faithful character of $E$ and $E$
is not abelian, we have that $\chi(1)>1$. Thus $\chi_Z$ is a 
reducible character.
Since $\chi_E\in \Irr(E)$, $E/Z$ is an abelian chief factor,
$\chi_Z$ is a reducible character  and 
$Z={\bf Z}(G)$, we have that $E/Z$ is a fully ramified section
with respect to $\chi_E$ and $\lambda$, and  that 
$\chi(1)^2=|E:Z|$.
Thus
\begin{equation}\label{dividingintwo}
(\chi\overline{\chi})_E= 1_Z^G.
\end{equation}
Since $G$ acts on $(E/{\bf Z}(G))^{\#}$ 
with 2 orbits, $G$ acts on $\Irr(E \modu Z)^{\#}$ with 
2 orbits.

 By 
Clifford Theory we have that for any character $\alpha\in \Irr(G)$,
the irreducible constituents of $\alpha_E$ form a $G$-orbit.
Let $\alpha_1\in \Irr(G)^{\#}$ be a character such that 
$[\alpha_1, \chi \overline{\chi}]=m_1\neq 0$. 
By Lemma \ref{center} we have that 
$\Ker(\alpha_1)\geq Z$. 
Thus the irreducible constituents
of $(\alpha_1)_E$ lie in  $\Irr(E \modu Z)^{\#}$.  
Since $G$ acts on  $\Irr(E \modu Z)^{\#}$ with 2 orbits and 
the irreducible constituents of $(\alpha_1)_E$ form a 
$G$-orbit,  there exists
 $\delta \in \Irr(E\modu Z)$ such that
 $[\delta, (\alpha_1)_E]=0$. Then there
exists a character $\alpha_2\in \Irr(G)$ such that 
$[(\alpha_2)_E, \delta] \neq 0$ and 
$[\alpha_2, \chi\overline{\chi}]=m_2\neq 0$. Thus the irreducible 
constituents of $(\alpha_1)_E$ and $(\alpha_2)_E$
form  two distinct $G$-orbits of $\Irr (E\modu Z)^{\#}$.
Since $G$ acts on  $\Irr(E \modu Z)^{\#}$ with 2 orbits
we conclude that 
 those are the only distinct $G$-orbits of $\Irr (E\modu Z)^{\#}$. 

Since $E/Z$ is an abelian chief
factor of $G$, by \eqref{dividingintwo}
we have that the irreducible constituents of $\chi\overline{\chi}$
appear with multiplicity 1. Since 
the irreducible 
constituents of $(\alpha_1)_E$ and $(\alpha_2)_E$
form  two distinct $G$-orbits of $\Irr (E\modu Z)^{\#}$ and the 
irreducible constituents of $\chi\overline{\chi}$ appear with multiplicity
1, we have that $\alpha_1$ and $\alpha_2$ are all
 the non-principal 
irreducible constituents of $\chi\overline{\chi}$ and they appear
with multiplicity 1. Thus
$$\chi\overline{\chi}=1+\alpha_1+\alpha_2.$$
\end{proof}

\begin{proposition}\label{examplekernelsequal}
There exist a group $G$ and a character $\chi \in \Irr(G)$ 
such that Hypothesis \ref{hypot}  holds and 
	$$\Ker(\alpha_1)=\Ker(\alpha_2).$$
Thus there exist a group $G$ and a character $\chi \in\Irr(G)$ 
satisfying 
the hypotheses of
Proposition \ref{kernelsequal}. 
\end{proposition}
\begin{proof}
Let $E$ be an extraspecial group of order $27$ and exponent $3$.
The  2-Sylow subgroup of 
$\SL(2,3)$ is isomorphic to the group ${\bf Q}_8$
by Lemma \ref{symplectic2}. Thus 
$\SL(2,3)$ contains a cyclic group of order 4.  Let $H$ be
a cyclic group of order 4. By Lemma \ref{symplectic2} (i)
and Lemma \ref{splitting}, we may assume that $H$ is a subgroup
of $\Aut(E\modu {\bf Z}(E))$. We can check that $H$ acts on 
$(E/Z)^{\#}$ with 2 orbits of the same size. 

Set $G= E\rtimes H$. Observe that $E/{\bf Z}(E)$ is an abelian 
chief factor of $G$. Also observe that $G/E$ acts faithfully
on $E/{\bf Z}(E)$, and with 2 orbits on $(E/{\bf Z}(E))^{\#}$.
Let $\Lambda\in \Irr(E)$ be a character of degree 3.
Observe that $\Lambda$ is $G$-invariant.
Since $(|H|, |E|)=1$ and $\Lambda$ is $G$-invariant, we have that
$\Lambda$ extends to $G$. Let $\chi \in \Irr(E)$ be an extension
of $\Lambda$. Observe that the hypotheses of Lemma \ref{basicequalkernel}
hold for $G$ and $\chi$. Thus the result follows from Lemma 
\ref{basicequalkernel}.
\end{proof}
\begin{proposition}\label{examplesabeliancase}
Let $p$ be a prime number and $m>0$ an
 integer.
Then there exist a solvable group $G$ and a  
character $\chi \in \Irr(G)$ such that $\chi(1)= p^m$ and
 the hypotheses of Proposition \ref{abelian} hold. Also, there exist
 a supersolvable group $G$ and a character $\chi\in \Irr(G)$ such 
 that $\chi(1)=p$ and the hypotheses of Proposition \ref{abelian} hold. 
\end{proposition}

\begin{proof}
	Let $\GF(q)$ be a finite field with $q$ elements, where
$q=p^m$. 
Denote by $F^{*}$ the group of units of  $\GF(q)$.
Also denote by $F$ the additive group of  $\GF(q)$.
 Observe that $F^{*}$ acts on $F$ by multiplication.

	Define $M$ to be the semi-direct product 
$F \rtimes F^{*}$ 
of $F$ by  $F^{*}$. Thus, for any $f_1, g_1 \in F$ and
$f_2, g_2 \in F^{*}$, we have
$$(g_1 \rtimes g_2)(f_1 \rtimes f_2)= (g_1 +g_2 f_1) \rtimes g_2 f_2.$$
Observe that if $m=1$ the group $M$ is a cyclic by cyclic group and thus
is supersolvable. 

Set $X=\GF(q)$. For any $f_1 \rtimes f_2 \in M$ and $x \in X$, define 
	$$(f_1 \rtimes f_2)(x) = f_1 + f_2 x.$$
\noindent We can check that this defines an action of $M$ on 
the set $X$.
We can also check that $0\in X$ is fixed by $F^{*}$.
 
Let $C$ be a cyclic group of order $p$. Set
$$K= C^{X} =  \mbox{ the set of all functions from $X$ to $C$}.$$
Observe that $K$ is just an elementary abelian group of order 
$p^{|X|}$. We can check 
 that $$(fm)(x)= f(mx),$$ for any $m \in M$ and $x \in X$, defines
an action of $M$ on $K$. 
	 
	Let $G$ be the wreath product of $C$ and $M$ relative to $X$,
i.e. $G = K \rtimes M$. Observe that if $m=1$, then the group $G$
is a supersolvable group. 

Let $\lambda \in \Irr(C)^{\#}$ be a character. 
Choose $\theta_0 \in \Irr(K)$ so that we have, for all $k \in K$,
$$\theta_0(k) =\lambda (k(0)), \mbox { where } 0=0_F \in X.$$ 
Since $0=0_F\in X$ is fixed by $F^{*}$, the stabilizer of 
$\theta_0$ in $G$ is $G_{\theta_0}= K\rtimes F^{*}$. We can check that
$\theta_0$ extends 
to $G_{\theta_0}$. Let ${\theta_0}^e \in \Irr( G_{\theta_0})$ be an
extension of $\theta_0$. Set $(\theta_0^e)^G= \chi$.
By Clifford theory we have that $\chi\in \Irr(G)$. 
Observe that  $\chi(1)=q$  since $G_{\theta_0}= F^{*} K$ 
and $\theta_0$ extends 
to $G_{\theta_0}$. 
Since ${\theta_0}^e\in \Irr(G_{\theta_0})$
 is an extension of $\theta_0$ and $(\theta_0^e)^G= \chi$,
we have that $[\chi_K, \theta_0]=1$. We can check that 
\begin{equation*}\label{ejemplo}
\chi_K= \sum_{x \in X}\theta_x
\end{equation*}
\noindent where $\theta_x(k) = \lambda (k(x))$ for all $k \in K$.
Therefore 
\begin{equation*}\label{ejemplosuma}
(\chi\overline{\chi})_K= \sum_{x, y \in X}\theta_x \overline{\theta_y}.
\end{equation*}
Since
$G_{\theta_0}= K\rtimes F^{*}$, $G= (K \rtimes F^{*})(K \rtimes F)$
  and
$(K\rtimes F^{*})\cap (K\rtimes F)= K$, we have that 
$\chi_ {K \rtimes F}= (\theta_0)^{K\rtimes F}\in \Irr(K\rtimes F)$.
Thus $\chi\overline{\chi}(g) =0$ for all $g \in (K\rtimes F) \setminus K$.
We can check that for any 
$\gamma \in \Irr( (K\rtimes F )\modu K)$,
\begin{equation}\label{gradoalpha1}
[\gamma,(\chi\overline{\chi})_{ K\rtimes F}]\neq 0.
\end{equation}
We can check that $\alpha_1=\gamma^G \in \Irr(G)$. Since 
$|G: K\rtimes F|=|F^{*}|=q-1$, we have that
$\alpha_1(1)=q-1$ . Observe  that
\eqref{gradoalpha1} implies
\begin{equation}\label{alpha1bajochi}
 [\alpha_1, \chi\overline{\chi}]\neq 0.
\end{equation}
Define $\delta= \theta_0 \overline{\theta_1}$.
Observe that $\theta_x (1)=1$ for all $x \in X$. Therefore 
$\delta \in \Irr(K)$. Since $0$ and $1$ are not fixed simultaneously
by any non-trivial element of $M$, we have that
$G_{\delta}=K$. Thus $\delta^G= \alpha_2 \in \Irr(G)$.
Observe that 
\begin{equation}\label{alpha2bajochi}
[\alpha_2, \chi\overline{\chi}]\neq 0
\end{equation}
\noindent since $[\delta, (\chi\overline{\chi})_K]\neq 0$.
Observe that $\alpha_2(1)=q(q-1)=|G: K|$.

 Observe that $\chi(1)=q$, $\alpha_1(1)=q-1$, $\alpha_2(1)=q(q-1)$,
 \eqref{alpha1bajochi} and \eqref{alpha2bajochi} imply
$$\chi\overline{\chi}= 1_G +\alpha_1 + \alpha_2.$$
\end{proof}

\begin{proposition}
There exist a solvable group $G$ and a faithful 
character $\chi\in \Irr(G)$  such
that $\chi(1)=10$ and 
$$\chi\overline{\chi}= 1_G + \alpha_1 + 4 \alpha_2,$$
\noindent where $\alpha_1, \alpha_2 \in \Irr(G)$, $\alpha_1(1) = 3$
and 
 $\alpha_2(1)=24$. Thus there exist a group $G$ and a character
 $\chi\in \Irr(G)$ satisfying  the conditions in
case (i) of Proposition \ref{knotabelian}.
\end{proposition}

\begin{proof}
Let $E$ be an extra-special group of order $5^3$ and exponent
$5$. Then 
$E/{\bf Z}(E)$ is an elementary abelian group of order
25.
By Lemma \ref{symplectic2} (iv) there exists 
a subgroup $H$ of  $\SL(2,5)$  isomorphic
to $\SL(2,3)$.  Assume $H$ 
acts on the vector space $E/{\bf Z}(E)$ in the natural way.
 We can check
that in fact $H$ acts on $(E/{\bf Z }(E))^{\#}$ regularly
and transitively.
By Lemma \ref{symplectic2} (i), we have that 
 $\Spl(2,5) \cong \SL(2,5)$. By Lemma \ref{splitting}, 
there exists a exact sequence
\begin{equation*}
0 \rightarrow \Inn (E) \rightarrow \Aut(E\modu {\bf Z}(E))\rightarrow
 \Spl(2,5)\rightarrow 0.
\end{equation*}
Since $\Spl(2,5) \cong \SL(2,5)$, we may assume that  $H$ acts on $E$.
Define  $G= E\rtimes H$ to be the semi-direct product of $E$ by $H$. 
Then ${\bf Z}(G)={\bf Z}(E)$. 

Let $\Lambda \in \Irr(E)$ be a non-linear character. Observe that
$\Lambda$ is $G$-invariant. Since $H$ acts 
on $E$ with coprime
 action, 
 $\Lambda$ extends to $G$, i.e, there exists $\phi \in \Irr(G)$
such that $\phi_E=\Lambda$. Let $\gamma \in \Irr(G\modu E)$ be
 a character of degree 2. Set $\chi= \phi \gamma$.
Then $\chi \in \Irr(G)$ and $\chi(1)=10$.
     
	Observe that $\chi(e)=0$ for any $e \in E \setminus {\bf Z}(E)$.

	Let $\delta \in (\Lin( E \modu {\bf Z}))^{\#}$.
 Since $H$ acts on $(E/{\bf Z}(E))^{\#}$ regularly and 
 transitively, $\delta$ induces irreducibly. Let $\alpha_2= \delta^G \in \Irr(G)$.
Observe that  
$\alpha_2(1)= 24$. Also  $\alpha_2 (g)=0$ if $ g\in G\setminus E$.
Thus 
\begin{equation*}
\begin{split}
[\chi\overline{\chi}, \alpha_2] & =
\frac{1}{|G|} \sum_{g \in G} (\chi\overline{\chi})(g)\
 \overline{\alpha_2(g)}\\
&= \frac{1}{|G|}[ \sum_{g \in G\setminus E}
(\chi\overline{\chi})(g)\  \overline{\alpha_2(g)} +
\sum_{g \in E \setminus {\bf Z}(G) }(\chi\overline{\chi})(g)\
 \overline{\alpha_2(g)}
+ \sum_{g \in {\bf Z}(G)}(\chi\overline{\chi})(g)\ \overline{\alpha_2(g)}]\\
&=\frac{1}{|G|} \sum_{g \in {\bf Z}(G)}(\chi\overline{\chi})(g)\  \overline{\alpha_2(g)}
=4.
\end{split}
\end{equation*}

We can check that $\gamma\overline{\gamma}= 1+ \alpha_1$,
where $\alpha_1(1)=3$ and $\alpha_1 \in \Irr(G \modu E)$.
Thus
\begin{equation*}
\begin{split}
[\chi\overline{\chi}, \alpha_1]& = [\chi, \chi\alpha_1]\\
				&= [ \phi\gamma, \phi \gamma \alpha_1]\\
		&= [ \phi\gamma\overline{\gamma}, \phi\alpha_1]=1
\end{split}
\end{equation*}
\noindent where the last inequality is since 
$[\gamma\overline{\gamma}, \alpha_1]=1$.

	Thus $\chi{\overline{\chi}}= 1_G + \alpha_1 + 4\alpha_2$, 
since 
$(\chi\overline{\chi})(1)= 10^2= 1 + 3+ 96 = 1+\alpha_1(1) + \alpha_2(1)$.
\end{proof}

\begin{proposition}
There exist a solvable group $G$ and a faithful 
character $\chi\in \Irr(G)$ such
that
$$\chi\overline{\chi}= 1_G + \alpha_1 + 4 \alpha_2$$
\noindent where $\alpha_1, \alpha_2 \in \Irr(G)$, $\chi(1)=14$,  
$\alpha_1(1) = 3$
and $\alpha_2(1)=48$. Thus there exist a group $G$ and a character 
$\chi \in \Irr(G)$ satisfying the conditions in case (ii) of
 Proposition \ref{knotabelian}.
\end{proposition}
\begin{proof}
Let $E$ be an extra-special group of order $7^3$ and exponent 7.
Observe that $E/{\bf Z}(E)$ is an elementary abelian 7-group. 
By Lemma \ref{sl27}, we may regard
$\widetilde{\GL}(2,3)$ as  a subgroup of $\SL(2,7)$. 
 By Lemma \ref{symplectic2} and Lemma \ref{splitting}, 
we may assume that $\widetilde{\GL}(2,3)$ acts on $E$.
Define $G$ to be 
 the semi-direct product  $E\rtimes  \widetilde{\GL}(2,3)$ of $E$
by $\widetilde{\GL}(2,3)$. We can check as before that there exists
a faithful 
character $\chi\in \Irr(G)$ satisfying the conditions of the proposition.
\end{proof}
\end{subsection}
\begin{subsection}{Dihedral and nilpotent groups}

Using the character table of the dihedral group of 
order $2^m$ with $m>3$,
we can check that dihedral groups have characters $\chi$ that satisfy 
Hypothesis \ref{hypot}.
Observe that therefore the index $|G:{\bf Z}(G)|$ is not 
going to be bounded as it
was in the case where $\chi\overline{\chi}=1_G +m\alpha$.
We are going to see now that this is, in some sense,
 the general case 
whenever we 
are working with nilpotent groups.
\begin{theorem}\label{nilpotent}

	Assume Hypothesis \ref{hypot}. Assume also that $G$ is nilpotent.
Then 
$$m_1=m_2=1, \  \ \chi(1)=2$$
\noindent and $G/ {\bf Z}(G)$ is isomorphic to a dihedral 2-group.
\end{theorem}
\begin{proof} By Lemma \ref{kernelcenter}, either  
$\Ker(\alpha_1)=Z$ or $\Ker(\alpha_2)=Z$.
Since the chief factors of a nilpotent group are 
cyclic of prime order, we have by Propositions \ref{kernelsequal},
\ref{knotabelian} and \ref{abelian} that either 
$\Ker(\alpha_1)\neq Z$ and $\Ker(\alpha_1)$ is abelian,
 or $\Ker(\alpha_2)\neq Z$ and  $\Ker(\alpha_2)$ is abelian.
Without 
lost of generality we may  assume
that $\Ker(\alpha_1)\neq Z$ is an abelian group  and $\Ker(\alpha_2)=Z$.  

We use now the notation of Lemma \ref{abelian}.
Then $|G:G_{\theta}|= |L:K|$ is a prime number. Since $G$ is nilpotent,  
$G_{\theta}$ is a normal subgroup of $G$. 
Observe that $K \leq  G_{\theta}$.
If $G_{\theta}\neq K$, then
 $\chi_{ G_{\theta}}\in \Irr(G_{\theta})$ by 
Lemma \ref{lemma1}. But this 
can not be since $\chi= \chi_{\theta}^G$ and $G_{\theta} < G$. 
Thus $G_{\theta}=K$ and so 
$G=L$.

	Since $K=\Ker(\alpha_1)$ is abelian and $|G:K|=p$ is a prime number, 
 $\alpha_1(1)=1 $ and 
either $\alpha_2(1)=1$ 
or $\alpha_2(1)=p$. Assume that  $\alpha_2(1)=1$. 
Since $\Ker(\alpha_2)=Z$, the group  $G/Z$ is cyclic. Since $Z={\bf Z}(G)$ and
$G/Z$ is cyclic, it follows that $G$ is an abelian group. Thus
$\chi(1)=1$.
But then  $\chi$ can not satisfy \eqref{basic}. So $\alpha_2(1)\neq 1$.
We conclude that $\alpha_2(1)=p$.

	Since $\chi$ and $\alpha_2$ are induced from characters of $K$,
 we have
$\chi(g)=\alpha_2(g)=0$ for all $g \in G \setminus  K$. Thus
$\alpha_1(g)=-1/a_1=-1$ for all $g \in G \setminus K$. Since $G/K$ 
is a cyclic
group of order $p$, $\alpha_1(g)=-1$ for all $g \in G \setminus K$ and 
$\alpha_1(1)=1$,
we must have that $p=2$.
Since $\chi(1)=p$, $\alpha_2(1)=p$ and 
$\chi(1)^2 = 1 + \alpha_1(1)+ m_2 \alpha_2(1)$, we conclude
that $m_2=1$.
	
	We want to see now that $G/Z$ is isomorphic to a dihedral 
2-group.
Let $\chi_K= \theta_1+ \theta_2$ and 
$(\alpha_2)_K= \gamma_1 + \gamma_2$, where
 $\theta_1, \theta_2,  \gamma_1, \gamma_2\in \Irr(K)$, 
$\theta_1 \neq \theta_2$ and $\gamma_1 \neq \gamma_2$.
This holds
  because $\chi$ and $\alpha_2$ are induced from characters
of the abelian normal subgroup $K$ and $|G:K|=2$.
Then by \eqref{basic} we have that 
\begin{eqnarray*}  \chi_K \overline{\chi_K} & = & 
(\theta_1+\theta_2)\overline{(\theta_1+\theta_2)} \\
            & = & 
2\, 1_K + \theta_1\overline{\theta_2} + \overline{\theta_1}\theta_2 \\
                    & =& 
2\, 1_K + (\alpha_2)_K \\
         & =&
2 \, 1_K + \gamma_1 + \gamma_2.
\end{eqnarray*}
Choose $\gamma_1 = \theta_1\overline{\theta_2}$ and 
$\gamma_2 = \overline{\theta_1}{\theta_2}$.
 Observe that 
$\theta_1\overline{\theta_2}= \overline{\overline{\theta_1}\theta_2}$.
Thus $\gamma_1= \overline{\gamma_2}$.
If $\gamma_1$ is a  real character, then 
$\gamma_1= \overline{\gamma_1}= \gamma_2$. But $\gamma_1\neq \gamma_2$.
Therefore $\gamma_1$ and $\gamma_2$ are not real characters.

 Since $\gamma_1$ is a linear character,
it follows that $(\gamma_1)^{-1}= \overline{\gamma_1}= \gamma_2$. 
In particular $\Ker(\gamma_1) = \Ker(\gamma_2)$. Therefore 
 $Z= \Ker(\alpha_2)= \Ker(\gamma_1) \cap \Ker(\gamma_2) = \Ker(\gamma_1)$.
Since $\gamma_1(1)=1$ , the group 
$K/Z$ is  cyclic. Let $\bar{x}$ be a generator of 
$K/Z$. Let $y$ be a coset representative of $G\setminus K$. 
Since $(\theta)_K= \gamma_1 + \gamma_2$, $\gamma_1 \neq \gamma_2$ and 
$y$ is a coset representative of $G\setminus K$, by Clifford
theory it follows that $(\gamma_1)^y= \gamma_2$.
Since $\gamma_1^{y^2}=\gamma_1$, $\gamma_1^y=\gamma_2= (\gamma_1)^{-1}$
and
$\gamma_1 \in \Irr(K \modu Z)$ is
a faithful character, we have that
\begin{equation}\label{dihedralrelation}
	\bar{y}\bar{x}\bar{y}^{-1}=\bar{x}^{-1}.
\end{equation}
Thus $G/Z$ is either isomorphic to the dihedral group or to 
the generalized quaternion group. If $G/Z$ is isomorphic to the 
generalized quaternion group, then $G$ is a
non-split central extension of 
a cyclic group by a generalized quaternion group. But that can not
be since $H^2({\bf Q}_{2^n}, {\bf C})= 1$ for any $n\geq 3$ (see
example in page 301 of \cite{reiner}).
Thus $G/Z$ is isomorphic to the dihedral 2-group and the proof
is complete.

\end{proof}
\end{subsection}
\end{section}


\begin{thebibliography}{9}

\bibitem{edithd} E. Adan-Bante, Products of characters and derived length, J. Algebra
266 (2003), 305-319.
 
\bibitem{reiner} C. W. Curtis, I. Reiner, Methods of Representation
Theory Vol 1 with applications to finite groups and orders, 
John Wiley $\&$ Sons  Inc, New York 1990.

\bibitem{go} D. Gorenstein, Finite Groups. New York: Harper and Row 1968.

\bibitem{huppert} B. Huppert, Endliche Gruppen I,
 Berlin--Heidelberg--New York : Springer 1967.

\bibitem{isaacs} I. M. Isaacs, Character Theory of Finite Groups. 
New York-San
Francisco--London: Academic Press 1976.

\bibitem{keller} T. M. Keller, Orbit sizes and character degrees II,
J. reine angew. Math. 516 (1999), 27-114.

\bibitem{robinson}
 D. J. S.  Robinson,  A Course in the Theory of Groups,
Springer-Verlag, New York 1996.


\bibitem{winter} D. L.  Winter, The automorphism group of an 
extraspecial $p$-group, Rocky mountain J. of Math 2
(1972), 159-168.

\bibitem{wolf} O. Manz, T. R. Wolf,  Representations of Solvable Groups, 
Lecture Notes Series Vol 185, Cambridge University Press 1993. 




\end{thebibliography}
\end{document}